\documentclass{article}

\usepackage[utf8]{inputenc}
\usepackage{amsmath,amssymb,amsfonts,amsthm}
\usepackage{hyperref}
\usepackage{color}
\usepackage{graphicx}
\usepackage{epstopdf}
\usepackage{tikz}
\usepackage{enumerate}

\usepackage{algorithm}
\usepackage{algorithmic}
\usepackage{array}
\usepackage{eqparbox}

\usepackage{etoolbox}  % patch def of algorithmic environment
\makeatletter
\patchcmd{\algorithmic}{\addtolength{\ALC@tlm}{\leftmargin} }{\addtolength{\ALC@tlm}{\leftmargin}}{}{}
\makeatother

\newcommand{\R}{\mathbb{R}}

\newcommand{\first}[1]{\mathsf{first}(#1)}

\newtheorem{theorem}{Theorem}[section]
\newtheorem{lemma}{Lemma}[section]
\newtheorem{definition}{Definition}[section]
\newtheorem{corollary}{Corollary}[section]
\newtheorem{proposition}{Proposition}[section]

\newtheorem{remark}{Remark}[section]

\newtheorem{assumption}{Assumption}[section]

\begin{document}

\begin{center}
{\Large \textbf{A new algorithm for the $^K$DMDGP subclass of Distance Geometry problems\footnote{This is an expanded and revised version of the abstract submitted to CTW2020.} }}
\end{center}

\begin{flushleft}
Douglas S. Gon\c{c}alves

\emph{CFM, Federal University of Santa Catarina}\newline 
\emph{88.040-900, Florian\'{o}polis - SC, Brazil} \newline
\texttt{douglas@mtm.ufsc.br} \\
\ \\
Carlile Lavor

\emph{IMECC, University of Campinas\newline
13081-970, Campinas - SP, Brazil} \newline
\texttt{clavor@ime.unicamp.br}\newline
\ \\
Leo Liberti

\emph{LIX CNRS \'{E}cole Polytechnique, Institut Polytechnique de Paris, \newline 
91128, Palaiseau, France} \newline 
\texttt{liberti@lix.polytechnique.fr}\newline
\ \\
Michael Souza

\emph{Federal University of Cear\'{a}\newline
60440-900, Fortaleza - CE, Brazil} \newline
\texttt{michael@ufc.br}
\end{flushleft}

\begin{abstract}
The fundamental inverse problem in distance geometry is the one of finding positions from interpoint distances. 
The Discretizable Molecular Distance Geometry Problem (DMDGP) is a subclass of
the Distance Geometry Problem (DGP) whose search space can be discretized
and represented by a binary tree, which can be explored by a
Branch-and-Prune (BP) algorithm. It turns out that this combinatorial search
space possesses many interesting symmetry properties that were studied in the
last decade. In this paper, we present a new algorithm for this subclass of the DGP, 
% for solving the DMDGP, which leverages symmetry information in this class of the DGP. 
which exploits DMDGP symmetries more effectively than its predecessors. 
Computational results show that the speedup, with respect to the
classic BP algorithm, is considerable for sparse DMDGP instances 
related to protein conformation.
\end{abstract}

\section{Introduction}\label{sec:intro}
Given a simple, undirected, weighted graph $G=(V,E,d)$, with weight function $d:E\rightarrow
(0,\infty )$ and an integer $K>0$, the \textit{Distance Geometry Problem}
(DGP) consists in finding a \textit{realization} $x:V\rightarrow \mathbb{R}^{K}$ such that, 
\begin{equation}
\forall \{u,v\}\in E,\text{ }||x_{u}-x_{v}||=d_{uv},  \label{eq:dgp}
\end{equation}%
where $||\cdot ||$ denotes the Euclidean norm, $x_{v}:=x(v), \forall v \in V$ 
and $d_{uv}:=d(\{u,v\}), \forall \{u,v\} \in E$. 
Each equation in \eqref{eq:dgp} is called a distance constraint. 
We say that a realization $x$ satisfies $d_{uv}$ if the corresponding distance constraint is verified. 
A realization satisfying all distance constraints in (\ref{eq:dgp}) is called a \textit{valid realization}. 
We shall call a pair $(G,K)$ a DGP instance. 

There are many applications of Distance Geometry, mainly related to $K\in \{1,2,3\}$ \cite{billinge16, billinge18, mucherino2013}. 
An application to Data Science can be found in \cite{liberti17}, and a very recent survey on this subject is \cite{liberti20}. 
An important class of the DGP arises in the context of 3D protein structure
calculations ($K=3$), with distance information provided by Nuclear Magnetic
Resonance (NMR) experiments \cite{crippen88, malliavin19, wutrich89}. 

Existence and uniqueness of DGP solutions, among other
theoretical aspects of the problem, are discussed in \cite{siam2014}.
Henceforth, we will consider that the DGP admits a solution.

\begin{assumption}\label{a:feas}
The solution set of \eqref{eq:dgp} is non-empty. 
\end{assumption}

%Although the DGP is NP-hard \cite{saxe}, using vertex orders from NMR data
%and properties of protein geometry \cite{cassioli2015, lavor19a, lavor2019}, it
%is possible to solve big instances by a combinatorial approach \cite{carvalho08}.

The DGP is naturally cast as a search in continuous space. 
Depending on the graph structure, however, combinatorial search algorithms can be defined, 
notably via the identification of appropriate vertex orders \cite{cassioli2015, lavor19a, lavor2019}. 
Although DGP is NP-hard \cite{saxe}, these combinatorial approaches 
allowed to show that it is Fixed Parameter Tractable (FPT) on certain graph structures, 
as those arising in protein conformation \cite{liberti2013}.

The aforementioned vertex orders define a DGP subclass, called the \textit{Discretizable
Molecular Distance Geometry Problem} (DMDGP) \cite{coap2012, lavor12},
formally given as follows.

\begin{definition}\label{def:dmdgp}
A DGP instance $(G,K)$ is a $^{K}$DMDGP if there is a vertex order $v_{1},...,v_{n}\in V$, such that
\begin{enumerate}
%\item[0.] \tred{The corresponding DGP with the subset of distances $d_{ij}$, such that $|i-j|\le K$, admits a solution.}

%\item[1.] $G[\{v_{i},...,v_{i+K}\}]$ is a clique, for $i=1,\dots,n-K$; 
\item[1.] $G[\{v_{1},...,v_{K}\}]$ is a clique; 

\item[2.] %$CM(v_{i-1},...,v_{i-K}) \ne 0$, for every $i>K$.
\begin{enumerate}
\item For every $i>K$, $v_{i}$ is adjacent to $%
v_{i-1},,...,v_{i-K}$,
\item $CM(v_{i-1},...,v_{i-K}) \ne 0$.
\end{enumerate}
\end{enumerate}
\end{definition}

In the above definition, $G[\cdot ]$ denotes the induced subgraph and $CM(v_{i-1},...,v_{i-K})$ is the Cayley-Menger
determinant of  $v_{i-1},...,v_{i-K}$ \cite[Sec.~2]{siam2014}. Its squared value is proportional to the $(K-1)$-volume of a 
realization $x_{i-1},\dots,x_{i-K}$ for $v_{i-1},...,v_{i-K}$. Condition $CM(v_{i-1},...,v_{i-K}) \ne 0$ means 
that the points $x_{i-1},\dots,x_{i-K}$ span an affine subspace of dimension $K-1$.

Although Definition~\ref{def:dmdgp} applies to any dimension $K$, 
therefore covering other applications rather than molecular conformation where $K=3$, 
the term ``molecular'' is commonly kept in the related literature \cite{coap2012,siam2014,cassioli2015}, regardless of the dimension, 
to enforce the property that the adjacent predecessors of $v_i$ are {\it contiguous} 
(the term  ``contiguous $K$-lateration order'' to mean $^K$DMDGP is used in \cite{cassioli2015}), 
a desirable property when ordering atoms of a protein \cite{coap2012,lavor19a}. 

%Given a vertex order $v_1, \dots, v_n$, 
%sometimes we abuse notation and use ``vertex $i$'' for $v_i \in V$ and ``edge $\{i,j\}$'' for $\{v_i,v_j\} \in E$.  
%Although Definition~\ref{def:dmdgp} considers an arbitrary dimension $K$, we will focus
%on the application related to protein conformation, where $K=3$, and denote $%
%^{3}$DMDGP simply by DMDGP. We also remark that to better illustrate the
%associated theory, some figures are given for $K=2$.
When the dimension $K$ is clear from the context, we shall simply use DMDGP rather than $^K$DMDGP.
Moreover, without loss of generality, whenever we denote an edge by $\{v_i,v_j\} \in E$, we will assume that $i < j$, 
i.e~$v_i$ precedes $v_j$ in the vertex order of Definition~\ref{def:dmdgp}. 

%For $K=3$, properties 1 and 2 imply that $G$ is composed by a chain of
%contiguous 4-cliques. 

Properties 1 and 2(a) of Definition~\ref{def:dmdgp} says that $G$ is composed by a chain of contiguous $(K+1)$-cliques.
Moreover, properties 1 and 2 allow us to turn the search space into a binary tree, in the following way. 

After fixing the positions for the first $K$ vertices, for each new vertex 
$v_{i}$, with $i>K$, property 2(a) ensures that the possible positions $x_{i}$
for vertex $v_{i}$ lie in the intersection of $K$ spheres centered at $x_{i-1},\dots,x_{i-K}$ 
with radii $d_{i-1,i},\dots,d_{i-K,i}$, respectively. 
Property 2(b) guarantees that there are at most two points, let us say $\{x_{i}^{+},x_{i}^{-}\}$, in such intersection \cite{maioli2017}. 
This spheres intersection can be computed in many different ways 
that we will not cover in this paper but are well studied in the literature \cite{alencar2019,maioli2017}.

\begin{remark}\label{rem:Klat}
The above process is known in the literature as $K$-lateration \cite{siam2014}. 
\end{remark}

Thus, following the vertex order, after fixing the first $K$ vertices,
each new vertex has at most 2 possible positions, which of course depend on
the position of its $K$ immediate adjacent predecessors, leading to a
binary tree of possible positions, where each path, from the root to a leaf
node, corresponds to a {\it possible realization} for the graph $G$. 

However, not all of these possible realizations (paths on the tree) are
valid, because $G$ may contain other edges $\{h,i\}$, with $|h - i| > K$, associated to
distance constraints that are not satisfied by such realizations. 
The edges given in Definition~\ref{def:dmdgp} are called \textit{discretization edges} and the
others, that may be (or not) available, are called \textit{pruning edges}. 

Henceforth, let us partition $E = E_D \cup E_P$, where $E_D$ is the set of discretization edges and $E_P$ the set of pruning edges.  
Clearly, we can also partition the equations in \eqref{eq:dgp} in discretization edge constraints and pruning edge constraints. 
We remark that, according to Definition~\ref{def:dmdgp}, $E_D = \{ \{i,j\} \in E \ | \ |i - j| \le K \}$ 
and therefore $E_P = \{ \{i,j\} \in E \ | \ |i - j| > K \}$.

The \textit{Branch-and-Prune} (BP) algorithm \cite{llm08}\
explores the DMDGP binary tree in a depth first manner and validates 
possible positions for vertices as soon as a pruning edge appears.   
A pseudo-code is given in Algorithm~\ref{algo:bp}. 

In Algorithm~\ref{algo:bp}, the phrase ``$x_i^+$ is feasible'' means that the equations 
$$
\forall h \ : \ h < i, \{v_h,v_i\} \in E_P,   \quad \| x_i^+ - x_h \| = d_{hi},
$$ 
are satisfied up to a certain tolerance. In Step~\ref{s:inter} positions $x_i^+$ and $x_i^-$ 
are computed via $K$-lateration. See \cite{siam2014,alencar2019,maioli2017} for details.
%The spheres intersection of Step~\ref{s:inter} can be computed in many different ways 
%that we will not cover in this paper but are well studied in the literature \cite{alencar2019,maioli2017}.

\begin{algorithm}[t]
\footnotesize
\begin{algorithmic}[1]
\STATE \textsf{BP}$(i,n,G,x)$ \COMMENT{($i>K$)}
\IF{($i > n$)} 
   \RETURN $x$ %\hfill \% {\it a solution is found}
\ELSE
%      \STATE Compute the intersection of spheres centered at $x_{i-K},\dots,x_{i-1}$, with radii $d_{i-K,i},\dots,d_{i-1,i}$ 
%      to obtain the candidate positions $\{x_i^+, x_i^{-}\}$. \label{s:inter}
      \STATE Find solutions $\{x_i^+, x_i^{-}\}$ for the system: $\| x_{\ell} - x_i \|^2 = d_{\ell,i}^2, \ell = i-K,\dots,i-1.$ \label{s:inter}
		\IF{$x_i^{+}$ is feasible}
		\STATE Set $x_i = x_i^{+}$ and call \textsf{BP}$(i+1,n,G,x)$. \COMMENT{1$^{st}$ candidate position}   
		\ENDIF
		\IF{$x_i^{-}$ is feasible}
		\STATE Set $x_i = x_i^{-}$ and call \textsf{BP}$(i+1,n,G,x)$. \COMMENT{2$^{nd}$ candidate position}   
	 \ENDIF
\ENDIF
\end{algorithmic}
\caption{BP}  \label{algo:bp}
\end{algorithm} 

Computational experiments in \cite{coap2012} showed that BP outperforms methods based on global continuation \cite{dgsol} and semidefinite programming \cite{kw2010} on instances of the DMDGP subclass, suggesting BP as the method of choice for this subclass of DGPs. 

In addition to the discretization of the DGP search space, the DMDGP order
also implies \textit{symmetry} properties of such discrete space \cite{powerof2-conf,gsi13,nucci13}. 
From the computational point of view, one of the most important of such properties, in the context of this paper, 
is that all DMDGP solutions can be determined from just one solution. %, found by any method. 
This property is related to the DMDGP symmetry vertices, which can be identified \textit{a priori}, based on the input graph (see next section). 
Once a first solution is found, the others can be obtained by partial reflections of the first, 
%based on these symmetry hyperplanes. 
based on symmetry hyperplanes associated to these vertices. 

%All of the previous works that exploit DMDGP symmetries consider that a
%DMDGP solution is given \cite{jbcb2012, liberti2013, jogo2018}. 
Previous works \cite{jbcb2012, liberti2013} exploited symmetry to reconstruct all valid realizations from the first one found and to prove that the BP algorithm is fixed-parameter tractable. Others \cite{pco12,jogo2018}, considered decomposition-based variants of BP 
which leverage DMDGP symmetry information. 
%For the first time, symmetry properties are used to find the ``first'' solution.

%In this work, we propose a new algorithm which exploits DMDGP symmetry in order to find the first valid realization more quickly.
In this work, we exploit DMDGP symmetry in order to find the first valid realization more quickly.
We handle  the DMDGP as a sequence of nested subproblems, each one defined 
%Although the proposed algorithm requires the DMDGP structure, it is intrinsically different from BP: 
%it handles the DMDGP as a sequence of subproblems, each one defined 
by a pruning edge $\{i,j\} \in E_P$. 
%and once a subproblem is solved, it is then considered as a rigid body. 
For each subproblem, we can exploit {\it any} realization $x$ (valid or not) for building the symmetry hyperplanes (which will define partial reflections).  
Once we have them, we apply compositions of such partial reflections {\it only} to $x_j$ to find its correct position. 
Only after finding the correct combination of partial reflections do we use it to obtain the positions of other vertices. 
After a subproblem is solved, the set of valid partial reflections is reduced and 
a single symmetry hyperplane is enough to handle positions $x_{i+K},\dots,x_j$ in the next subproblem.

In terms of the system of nonlinear equations \eqref{eq:dgp}, 
we solve a subset of equations and then gradually include new equations to this subset: 
the new equations are solved subject to the original equations in the subset. 
This process is repeated until all equations in \eqref{eq:dgp} are satisfied. 

These ideas lead to a new algorithm which deals with pruning edges, one-by-one, 
%in a specific order that makes the synchronization (merging solutions of subproblems) easier 
and takes advantage of a valid realization for already solved subproblems. 
Computational results illustrate the advantage of the new algorithm, compared with the classic BP.

This paper is organized as follows. Section~\ref{sec:sym} briefly reviews the main
results about DMDGP symmetries and Section~\ref{sec:sub} explains how they can be used
to solve a sequence of nested subproblems. The new algorithm, its correctness and implementation details are presented 
in Section~\ref{sec:algo}, and comparisons with the classic BP in protein-like instances
are given in Section~\ref{sec:results}. Concluding remarks are given in Section~\ref{sec:final}.

\section{DMDGP symmetries}\label{sec:sym} 
Before discussing the new algorithm, we shall present a theoretical background on DMDGP symmetries and 
recall some results from \cite{liberti2013,dam2014,siam2014}.  
%a corollary of the main result related to DMDGP symmetries, discussed below.

Given a realization $x$ satisfying \eqref{eq:dgp}, it is clear that 
there are uncountably many others, 
%which satisfy the same set of distances, that may be obtained 
which satisfy the same set of distances, and which can be obtained  
by translations, rotations or reflections of $x$ (because these transformations preserve all pairwise distances). 
Since the assumptions of Definition~\ref{def:dmdgp} ensure that the first $K$ vertices form a clique, 
a valid realization for $G[v_1,\dots,v_K]$ in $\R^K$ can be found by matrix decomposition methods \cite{dokmanic2015} 
or a sequence of spheres intersections \cite{alencar2019}, for example. 
Once the positions of these first $K$ vertices are fixed, the degrees of freedom of translations and rotations are removed. 

From here on, we say that two realizations are incongruent (modulo translations and rotations) 
if they are not translations, rotations or total reflections of each other.  
For technical reasons, we only allow the total reflections through the hyperplane defined by 
the positions of the first $K$ vertices in the vertex order 
(so two realizations, one of which is a reflection of the other through this hyperplane, will both be considered members of any set of ``incongruent realizations'').

\begin{definition}\label{def:hatX}
Let $\hat{X}$ be the set of all incongruent %(modulo translations and rotations)  
realizations satisfying distance constraints associated to discretization edges in $E_D$, i.e., $\|x_i-x_j\|=d_{ij}$ such that $|i - j| \le K$.  
A realization $x \in \hat{X}$ is called a possible realization. 
\end{definition}

As discussed in Section~\ref{sec:intro}, each $x \in \hat{X}$ corresponds to a path from the root to a leaf node 
in the binary tree of a DMDGP instance. Notice that $|\hat{X}| = 2^{|V| - K}$. 

\begin{definition}\label{def:X}
A realization $x \in \hat{X}$ is said to be valid if $x$ is a solution of \eqref{eq:dgp}. 
Let $X \subset \hat{X}$ denote the set of all incongruent %(modulo translations and rotations) 
valid realizations of a DMDGP instance.
\end{definition}

\begin{figure}
\centering
\begin{tikzpicture}
\draw [gray,dashed,domain=200:340] plot ({sqrt(5)*cos(\x)}, {sqrt(5)*sin(\x)}); 
\draw [black,thick,domain=200:340] plot ({sqrt(8)*cos(\x)}, {sqrt(8)*sin(\x)}); 
\draw [black,thin,domain=200:340] plot ({sqrt(10)*cos(\x)}, {sqrt(10)*sin(\x)}); 
\node (1) at (0,0) [fill, circle, minimum size=2pt, inner sep=2pt] {};
\node (2) at (0,-1) [fill, circle, minimum size=2pt, inner sep=2pt] {};
\node (31) at (-1,-2) [fill, circle, minimum size=2pt, inner sep=2pt] {};
\node (32) at (1,-2) [fill, circle, minimum size=2pt, inner sep=2pt] {};
\node (a) at (-2,-2) [fill, circle, minimum size=2pt, inner sep=2pt, label=below:{$a$}] {};
\node (d) at (2,-2) [fill, circle, minimum size=2pt, inner sep=2pt, label=right:{$d=R_{x'}^4(c) = R_x^3(a)$}] {};
\node (b) at (-1,-3) [fill, circle, minimum size=2pt, inner sep=2pt, label=below:{$b=R_x^4(a)$}] {};
\node (c) at (1,-3) [fill, circle, minimum size=2pt, inner sep=2pt, label=-60:{$c=R_x^3(b)=R_x^3(R_x^4(a))$}] {};
\draw (1)--(2)--(31)--(a);
%\draw[blue,thick] (1)--(a);
%\draw[blue!40!white,thick] (1)--(b);
\draw[dashed] (1)--(2)--(32)--(d);
\draw[dotted] (0,1) -- (0,-4);
\draw[dotted] (1,0) -- (-2,-3);
\draw[dotted] (-1,0) -- (2,-3);
\end{tikzpicture}
\caption{The leftmost path/realization $x$ is represented by a straight line whereas the rightmost $x'$ by a dashed line. 
All 4 possible positions for the fourth vertex (denoted by $a,b,c$ and $d$) can be generated by $x$ and its induced reflections $R_x^3$ and $R_x^4$. 
Illustration for $K=2$.}\label{fig:reflections}
\end{figure}
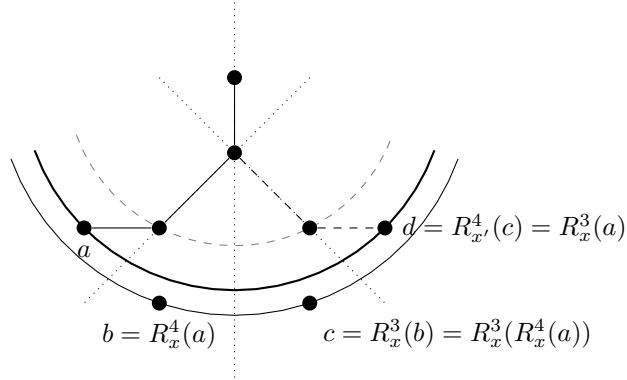

%Since the first DMDGP computational experiments \cite{coap2012},  it was noticed 
The computational experiments in \cite{coap2012} suggested
%that the number of valid realizations up to rotations and translations 
that $|X|$ is always a power of 2. 
A conjecture was formulated and quickly disproved using some instances constructed by hand,  
%This result remained unexplained for a while, until it was formally shown in \cite%
until the conjecture was shown to be true with probability one in \cite{dam2014}. 

Given $x \in \hat{X}$, for $i>K$, let $R_{x}^{i}(y)$ be the reflection of $y \in \R^K$ 
through the hyperplane defined by $x_{i-K},\dots,x_{i-1}$, with normal $p_i$:
$$
R_{x}^{i}(y) = (I - 2 p_i p_i^T)(y - x_{i-1}) + x_{i-1}, 
$$
assuming $\| p_i \| = 1$. 
%%% partial reflections %%%%%
Let us also define, for all $i>K$ and $x \in \hat{X}$, {\it partial reflection} operators:
\begin{equation}\label{eq:g}
g_i(x) = (x_1, x_2, \dots, x_{i-1}, R_x^i(x_i),R_x^i(x_{i+1}), \dots, R_x^i(x_n)). 
\end{equation}
\begin{remark}\label{rem:pr}
Some direct but useful properties of reflections and partial reflections are in order:
\begin{enumerate}
\item A reflection $R_x^i(y)$ preserves the distance from $y$ to any point in the hyperplane 
defined by $x_{i-K},\dots,x_{i-1}$. 
\item The pairwise distances for $R_x^i(x_i),R_x^i(x_{i+1}), \dots, R_x^i(x_n)$ are the same as those for $x_i,\dots,x_n$. 
As a consequence of this, and the fact that $R_x^i (x_{\ell}) = x_{\ell}$, for $\ell = i-K,\dots, i-1$, 
all pairwise distances for $x_{i-K},\dots,x_n$ from $x \in \hat{X}$ are preserved in $g_i(x)$. 
\item Partial reflections preserve distances related to discretization edges $E_D$, 
so that $g_i(x) \in \hat{X}$, for every $x \in \hat{X}$. 
\item All realizations in $\hat{X}$ can be generated from a single $x \in \hat{X}$ 
by the composition of partial reflection operators $g_i$ \cite[Sec.~3.3.8]{siam2014}.  
\end{enumerate}
\end{remark}

Let us now recall one of the main results about $^K$DMDGP symmetries.

\begin{theorem}[Theorem 3.2 in \cite{siam2014}]\label{thm:twos}
With probability 1, for all $j > K + i$, %$j > K$ and $i < j - K$ 
there is a set $H^{ij}$ of $2^{j-i-K}$ real positive values 
such that for each $x \in \hat{X}$, we have $\| x_j - x_i \| \in H^{ij}$. 
Furthermore, for all $x', x \in \hat{X}$ such that $x'\ne x$ and $x'_t = x_t$, for $t \le i+K-1$, 
$\| x_j - x_i \| = \| x_j' - x_i \|$ if and only if $x_j' = R_{x}^{i+K}(x_j)$.
%$\forall x \in \hat{X}$, 
%$\| x_j - x_i \| = \| R_{x}^{i+K}(x_j) - x_i \|$ and $\forall x' \in \hat{X}$, 
%if $x_j' \notin \{ x_j, R_{x}^{i+K} (x_j) \}$ then $\| x_j - x_i \| \ne \| x_j' - x_i' \|$.
\end{theorem}

In Theorem~\ref{thm:twos}, ``with probability 1'' means that the set of $^K$DMDGP instances 
for which the statements do not hold has Lebesgue measure zero in the set of all $^K$DMDGP instances \cite{dam2014}.
%%In Theorem~\ref{thm:twos}, ``with probability 1'' means: (i) that the intersection of 
%%spheres centered at $x_{i-K},\dots,x_{i-1}$ with radii $d_{i-K,i},\dots,d_{i-1,i}$ 
%%is composed by two distinct points $x_i^-$ and $x_i^+$, for all $i>K$, 
%%with probability one among the feasible instances of $^K$DMDGP
%%(because the set of feasible instances violating this has Lebesgue measure zero on the set of YES $^K$DMDGP instances);  
%%(ii) if $\{i,j\}$ is a pruning edge, the hyperplane defined by $x_i,\dots,x_{i+K-1}$ does not contain consecutive 
%%points $x_{\ell},\dots,x_{\ell+K-1}$, for $\ell = i+K, \dots, j-K+1$, an event which occurs with probability zero. 

The first part of Theorem~\ref{thm:twos} says that, for $j> K+i$,  %$j> K$ and $j - i > K$, 
the possible realizations $x \in \hat{X}$ yield a set of $2^{j-i-K}$ distinct values for $\| x_i - x_j \|$. 
Let $\hat{X}_{i+K-1}(x)$ be the subset of possible realizations $x' \in \hat{X}$ that agree with $x$ 
in the first $i+K-1$ positions. 
Given a possible realization $x$, each of these $2^{j-i-K}$ distinct values is associated to a pair of $2^{j-i-K+1}$ possible positions for $x_j$ 
from realizations in $\hat{X}_{i+K-1}(x)$ (see Figure~\ref{fig:reflections} where possible values for $\| x_1 - x_3\|$ and $\|x_1 - x_4\|$
are represented by the radii of gray and, respectively, black arcs centered at $x_1$). 

Since $j - i > K$, if the distance $d_{ij}$ is available, it must be a pruning distance. 
In view of Assumption~\ref{a:feas}, then $\| x_i - x_j \| = d_{ij}$ for some $x \in \hat{X}$. 
Let such a $x$ define the set $\hat{X}_{i+K-1}(x)$. 
Now, from the second part of Theorem~\ref{thm:twos}, we have that 
among the possible realizations $x' \in \hat{X}_{i+K-1}(x)$, 
only those such that $x_j' \in  \{ x_j, R_{x}^{i+K}(x_j)\}$ are feasible with respect to $d_{ij}$. 
If $v_j$ is the last vertex in the order, then only two realizations in $\hat{X}_{i+K-1}(x)$ are feasible. 

For every DMDGP solution, there is another one symmetric to the hyperplane 
defined by the positions of the first $K$ vertices. Moreover, as a consequence of Theorem~\ref{thm:twos}, 
the number of solutions doubles for every other {\it symmetry vertex} belonging to the following
set \cite{dam2014}:

\begin{equation}\label{eq:S}
S:=\{v_{\ell} \in V \ | \ \nexists \{v_i,v_j\}\in E\text{ with }i+K<\ell\leq j\}.
\end{equation}

The vertex $v_{K+1}$ is always in $S$, because the first $K$ vertices define
a symmetry hyperplane. The other symmetry hyperplanes are given by the positions of 
$v_{i-K},\dots,v_{i-1}$, if $v_{i}\in S$, for $i>K+1$. 
As mentioned in the Section~\ref{sec:intro}, $S$ can be computed before solving a $^K$DMDGP instance,  
which implies that the number of solutions is known \textit{a priori}, and given by $2^{|S|}$, with probability one. 

\begin{theorem}[Theorem~3.4 in \cite{liberti2013}]\label{thm:2S}
Let $(G,K)$ be a feasible $^K$DMDGP and $S$ its set of symmetry vertices. 
Then, with probability 1, $|X| = 2^{|S|}$. 
\end{theorem}

The $2^{|S|}$ valid realizations are incongruent modulo translations and rotations, 
meaning that they differ one from another only by partial reflections (or a total reflection through the first symmetry hyperplane, as explained above). 

It is important to notice from \eqref{eq:S} that the addition of new pruning edges in $E$ 
may reduce the number of elements (symmetry vertices) in $S$. 

A direct consequence of Theorem~\ref{thm:2S} is the following corollary.

\begin{corollary}\label{cor:twos}
Let $(G,K)$ be a feasible $^K$DMDGP instance where $|V(G)| = n>K$. 
If $\{v_1,v_n\} \in E$, then $(G,K)$ has only two incongruent solutions  
which are reflections of each other through the symmetry hyperplane defined by the position of the first $K$ vertices. 
\end{corollary}

\begin{proof}
If $\{v_1,v_n\} \in E$, then $S = \{v_{K+1}\}$, which implies that the number of solutions is $2^{|S|} = 2^{1} = 2$. 
If one of these solutions is $x$, then the other is %$x' = g_{K+1}(x)$. 
$x'$, the reflection of $x$ through the hyperplane defined by $x_1,\dots,x_K$.
\end{proof}

%\begin{corollary}\label{cor:1n}
%Let $(G,K)$ be a feasible $^K$DMDGP instance where $|V(G)| = n>K$ such that $\{v_1,v_n\} \in E$. 
%If $x \in \hat{X}$ is such that $\| x_1 - x_n \| = d_{1n}$, then for every $\{v_i, v_j \} \in E$ it holds $\| x_i - x_j \| = d_{ij}$. 
%\end{corollary}

%Furthermore, given one solution $x \in X$, all $2^{|S|} - 1$ others can be build from $x$ 
%by a sequence of partial reflections through hyperplanes corresponding to vertices in $S$. 

%Another useful result for solving a fundamental DMDGP subproblem is the following. 

A result that will be useful ahead is given in Proposition~\ref{prop:reflect} and illustrated in Figure~\ref{fig:reflections}. 

\begin{proposition}[Lemma 4.2 in \cite{liberti2013}]\label{prop:reflect}
Let $x \in \hat{X}$, $k > i + 1$ and $p_i \ne p_k$ be the normals to the hyperplanes 
defining $R_x^i(\cdot)$ and $R_x^k (\cdot)$. 
If $y \in \R^K$ is not in the hyperplanes containing the origin and normal to $p_i,p_k$,  
then $R_{g_i(x)}^k (R_x^i( y )) = R_x^i ( R_x^k ( y ) ) $.
\end{proposition}

Proposition~\ref{prop:reflect} tells us that compositions of partial reflections that depend on more 
than one realization (e.g~$x$ and $x' := g_k(x)$) can be described in terms of reflections based on a single realization.
For example, for $k>i$, we have 
\begin{eqnarray*}
(g_k \circ g_i)(x)  & = & g_k(g_i(x)) \\
\ & = & g_k(x_1,\dots, x_{i-1},R_x^i(x_i), \dots, R_x^i(x_n)) \\
\ & = & (x_1,\dots, x_{i-1},R_x^i(x_i), \dots, R_x^i(x_{k-1}), R_{x'}^k(R_x^i(x_k)),\dots, R_{x'}^k(R_x^i(x_n)) \\
\ & = & (x_1,\dots, x_{i-1},R_x^i(x_i), \dots, R_x^i(x_{k-1}), R_{x}^i(R_x^k(x_k)),\dots, R_{x}^i(R_x^k(x_n)),
\end{eqnarray*}
where the last equality follows from Proposition~\ref{prop:reflect}. 
%This idea will be important for solving DMDGP subproblems discussed in Section~\ref{sec:sub}. 
%%%%%%%%

Therefore, for a DMDGP, given $x \in \hat{X}$, problem \eqref{eq:dgp} can be cast 
as finding a binary vector $s \in \{0,1\}^{n-K}$, such that 
\begin{equation}\label{eq:xs}
x(s) := U(x,s) = g_{K+1}^{s_1} \circ \dots \circ g_{n}^{s_{n-K}} (x) 
\end{equation}
satisfies $|| x_i(s) - x_j(s) || = d_{ij}$, for all $\{i,j\} \in E$. 
Here, $g_i^{1}(\cdot) = g_i(\cdot)$ and $g_i^{0}(\cdot) = I(\cdot)$, where $I(x) = x$. 
In Section~\ref{sec:sub} we shall explain how to efficiently perform the search of this binary vector 
taking into account DMDGP symmetry information.  

To close this section, let us describe how to generate other valid realization $x(s') \in X$ 
from a given one $x(s) \in X$. 
Let $x(s)$ be a valid realization for $(G,K)$. 
The vertices in the set $S$ determine which components of the binary vector $s \in \{0,1\}^{n-K}$ from \eqref{eq:xs} 
are allowed to change in order to obtain another valid realization for $(G,K)$. 
In other words, the search space for the new $s' \in \{0,1\}^{n-K}$ is reduced to
\begin{equation}\label{eq:B}
s' \in B := \{ s' \in \{0,1\}^{n-K} \ | \ s'_{\ell} = s_{\ell} \text{ if } v_{K+\ell} \notin S \}. 
\end{equation}

\begin{lemma}\label{lem:NaoEstraga0}
Let $S \ne \varnothing$ and $x(s)$ be a valid realization for $(G,K)$.  
For every $s' \in B$, $x(s') \in X$. 
\end{lemma}
\begin{proof}
Since $x(s')$ from Eq.~\eqref{eq:xs} involves only partial reflections, 
in view of Property 3 in Remark~\ref{rem:pr}, $x(s') \in \hat{X}$, i.e~$\| x_i(s') - x_j(s') \| = d_{ij}, \forall \{i,j\} \in E_D$. 

It remains to show that $x(s')$ does not violate distance constraints associated to pruning edges $\{i,j\} \in E_P$. 
Since the reflections are applied to positions $x_{\ell}$ such that $\ell \ge K+1$, 
edges $\{i,j\} \in E$ with $i < j \le K$ are not affected. 
Thus, assume that $K+1 \le j \le n$.  
%If $u \le i$, then for $\ell = i+K+1,\dots,w$ there exists $\{u,w\}$ 
%such that $u+K+1 \le \ell \le w$, which implies that $v_{i+K+1},\dots,v_{w} \not\in S^{ij}$, 
%meaning that the first symmetry vertex $v_{\ell}$ in $S^{ij}$ is such that $\ell \ge w+1$. 
%Thus, according to \eqref{eq:g} and \eqref{eq:xs}, partial reflections are not applied to $x_{i+K+1},\dots,x_{w}$ and 
%$\| x_u - x_w \| = d_{uw}$ continues to be satisfied. 
We have that $v_{i+K+1},\dots,v_j \not\in S$, and from \eqref{eq:xs} and \eqref{eq:g}, 
positions $x_{\ell},\dots,x_{i+K+1},\dots,x_j$ are updated by reflections 
$R_{x}^{\ell}(x_{\ell}), \dots, R_{x}^{\ell}(x_{i+K+1}), \dots, R_{x}^{\ell}(x_j)$, for $K+1 \le \ell \le i+K$ such that $v_{\ell} \in S$. 
Since either $i \le \ell - 1$, i.e~$x_i$ is in the hyperplane associated to $v_{\ell}$, 
or $i \ge \ell$, i.e.~$x_i$ comes after this hyperplane, 
in view of Remark~\ref{rem:pr}, Property 2, these reflections are such that $\| x_i(s') - x_j(s') \| = d_{ij}$. 
\end{proof}

\section{Nested DMDGP subproblems}\label{sec:sub}
Given a DMDGP instance, properties 1 and 2 of Definition~\ref{def:dmdgp} give rise to a 
rich symmetric structure for the corresponding DGP problem, %which presents lots of symmetries, 
as discussed in Section~\ref{sec:sym}. 

On one hand, the absence of pruning edges turns the DMDGP into a trivial problem, 
because {\it any} path from the root to a leaf node of the search tree corresponds to a valid realization, 
i.e $\hat{X} = X$, 
and all other solutions can be built by partial reflections. 
On the other hand, one of the most challenging  DMDGP instances to solve with BP is the one where the only pruning edge is $\{v_1,v_n\}$. 
In that case, feasibility can only be verified at a leaf node, 
and for a standard depth-first search (DFS), 
it may represent a costly backtracking process until the first valid realization is found. 

Differently, given $x \in \hat{X}$, % from previous algorithms for DMDGP, such as BP which explores the search tree in a depth-first manner, 
the present proposal is to iteratively handle the pruning edge constraints following a given order $<$ on $E_P$. 

As mentioned in Section~\ref{sec:sym} (after Theorem~\ref{thm:2S}), 
each pruning edge $\{i,j\}$ may reduce the set of valid partial reflection operations 
that can be applied to realizations of the vertices $v_{i+K},\ldots,v_j$. 
Thus, by keeping track of valid partial reflections (or equivalently their corresponding symmetry vertices), 
it is possible to consistently modify a given realization satisfying a subset of distance constraints to also satisfy a new pruning edge constraint. 
This process is repeated until all distance constraints are satisfied. 

For this, we enumerate edges in $E_P$ as $e_1,e_2,\dots,e_m$, with $m = |E_P|$, 
and use $e_{k} < e_{\ell}$ to mean that edge $e_{k}$ precedes $e_{\ell}$ in this order. 
We define the set of pruning edges preceding edge $\{i,j\}$ by
\begin{equation}\label{def:Pij}
P^{ij} := \{ \{ u, w \}  = e' \in E_P \ | \ e' < e = \{i,j\} \}.
\end{equation}

Then, we define a sequence of subproblems spanned by $\{i,j\} \in E_P$ following the above pruning edge order. 

\begin{definition}\label{def:subprob}
Let $(G,K)$ be a feasible $^{K}DMDGP$ with $G=(V,E,d)$. 
Let $G^{ij} = (V,E^{ij},d_{|E^{ij}})$, where $E^{ij} = E_D \cup P^{ij} \cup \{i,j\}$, $\{i,j\} \in E_P$ and $d_{|E^{ij}}$ is the restriction of $d$ to $E^{ij}$. 
We say that $(G^{ij},K)$ is a $^{K}DMDGP$ subproblem of $(G,K)$ spanned by pruning edge $\{i,j\}$. 
\end{definition}

It is clear that $(G^{ij},K)$ is itself a DMDGP problem. Let us denote by $X(G^{ij})$ the solution set of $(G^{ij},K)$. 
%For two graphs $G=(V,E)$ and $H=(W,F)$, we use notation $G \subset H$ to mean $V \subset W$ and $E \subset F$. 

\begin{proposition}\label{prop:subset}
Let $G=(V,E,d)$ and $H=(V,F,\hat{d})$ such that $(G,K)$ and $(H,K)$ are feasible $^K$DMDGPs. 
If $E \subset F$ and $d(\{i,j\}) = \hat{d}(\{i,j\}), \forall \{i,j\} \in E$, then $X(G) \supset X(H)$. 
\end{proposition}

Let $(G^{uw},K)$ and $(G^{ij},K)$ be DMDGP subproblems spanned by edges $\{u,w\}$ and $\{i,j\}$, respectively, 
such that $\{u,w\} < \{i,j\}$. In view of Proposition~\ref{prop:subset}, we have $X(G^{uw}) \supset X(G^{ij})$. 

Moreover, in this sequence of DMDGP subproblems, each time a new pruning edge is included, e.g~$E^{ij} = E^{uw} \cup \{i,j\}$, 
the set of symmetry vertices (see Eq.~\eqref{eq:S}) for $(G^{uw},K)$ may be reduced. 
This motivates us to define the set of {\it necessary} symmetry vertices for subproblem $(G^{ij},K)$ as:
\begin{equation}\label{eq:Sij}
S^{ij} = \{ v_{\ell} \in \{v_{i+K+1},\dots,v_j\}  \ | \ \not\exists \{ u, w \} \in P^{ij}, u+K < \ell \le w \}.
\end{equation}

Let $x(s)$ be the current realization which is valid for $(G^{uw},K)$ and let $e_{k+1} = \{i,j\} > \{u,w\} = e_k$. 
The vertices in the set $S^{ij}$ determine which components of the binary vector $s \in \{0,1\}^{n-K}$ from Eq.~\eqref{eq:xs} 
are allowed to change in order to obtain a valid realization for $(G^{ij},K)$. 
In other words, the search space for the new $s' \in \{0,1\}^{n-K}$ is reduced to
\begin{equation}\label{eq:Bij}
s' \in B^{ij} := \{ s' \in \{0,1\}^{n-K} \ | \ s'_{\ell} = s_{\ell} \text{ if } v_{i+K+\ell} \notin S^{ij} \}. 
\end{equation}

\begin{lemma}\label{lem:NaoEstraga}
Let $S^{ij} \ne \varnothing$ and $e_{k+1} = \{i,j\} > \{u,w\} = e_k$. 
Let $x(s)$ be a valid realization for $(G^{uw},K)$.  
For every $s' \in B^{ij}$, $x(s') \in X(G^{uw})$. 
\end{lemma}
\begin{proof}
The proof is similar to the one of Lemma~\ref{lem:NaoEstraga0} and therefore is left in the Appendix.
\end{proof}

\begin{remark}\label{rem:key}
From Proposition~\ref{prop:subset} and Lemma~\ref{lem:NaoEstraga}, 
if $e_{k+1} = \{i,j\} > \{u,w\} = e_k$, then,  
given $x(s) \in X(G^{uw}) \supset X(G^{ij})$, to obtain $x(s') \in X(G^{ij})$ 
it suffices to find $s' \in B^{ij}$ such that $\| x_i(s') - x_j(s') \| = d_{ij}$. 
\end{remark}

%Since pruning edges in $P^{ij}$ are taken into account, 
%we may have $\{u,w\} \in P^{ij}$ such that $u+K < \ell \le w$, for $\ell=i+K+1,\dots,h \le j$,  
%which implies that $v_{i+K+1},\dots,v_h \not\in S^{ij}$. 

%We remark that \eqref{eq:Sij} is different from the set \eqref{eq:S} for subproblem $(G^{ij},K)$ with edge $\{i,j\}$ removed. 
%Let us denote by $\bar{S}^{ij}$ the set of symmetry vertices based on \eqref{eq:S} for $(G^{ij},K)$ without edge $\{i,j\}$. 
%Theorem~\ref{thm:2S} says that there are $2^{|\bar{S}^{ij}|}$ valid realizations for such problem, 
%but pruning edges in $P^{ij} \setminus E(G^{ij})$ further reduce the number of valid realizations 
%compatible with edges in $P^{ij}$ to $2^{|S^{ij}|}$. 

%Therefore, in order to preserve all pairwise distance constraints in already solved subproblems,  
%%and \tred{reduce the degrees of freedom (number of free variables)},  
%%the only components of $s$ allowed to take values in $\{0,1\}$ 
%the only components of $s$ that can be assigned the value 1 
%are those associated to %necessary symmetry vertices in $S^{ij}$. 
%elements of $S^{ij}$. 

%We now explain how we actually compute a realization for a subproblem given positions initialized as in Proposition~\ref{prop:initx}.  

%Nevertheless, real life DMDGP instances, e.g~the ones from protein conformation problems, 
%in fact are composed by many %nested or overlapped  
%DMDGP subproblems of this kind. 
Furthermore, in the following we show that there is a unique $s' \in B^{ij}$ 
satisfying such condition. For this, let us recall a simple fact that follows from Definition~\ref{def:dmdgp}.

\begin{proposition}\label{prop:subdmdgp}
If $(G,K)$ is a $^K$DMDGP instance, so is $(G[v_i,\dots,v_j],K)$, for $j>K+i$.%$j \ge i+K$.
\end{proposition}

Thus, given a $^K$DMDGP instance $(G,K)$, any subgraph induced by at least $K+1$ 
consecutive (w.r.t. the vertex order) vertices of $V(G)$ is a $^K$DMDGP itself. 
Proposition~\ref{prop:subdmdgp} implies that each $\{v_i,v_j\} \in E_P$  
%defines a DMDGP subproblem which has two solutions, which are reflections of each other 
defines a DMDGP instance based on the subgraph $G[v_i,\dots,v_j]$. 

\begin{proposition}\label{prop:two}
Any DMDGP instance $(G[v_i,\dots,v_j],K)$ spanned by $\{v_i,v_j\} \in E_P$ has only two solutions.
\end{proposition}

\begin{proof}
It follows from Proposition~\ref{prop:subdmdgp} and Corollary~\ref{cor:twos}.
\end{proof}

Proposition~\ref{prop:two} says that each DMDGP instance $(G[v_i,\dots,v_j],K)$
spanned by a pruning edge $\{i,j\}$ has only two solutions, 
which are reflections of each other through the hyperplane defined by $x_i, \dots, x_{i+K-1}$. 
These two solutions correspond to a particular configuration of the components $s'_{i+K},\dots,s'_{j}$. 
The only difference between the two is the first component $s'_{i+K}$. 
Since $v_{i+K} \not\in S^{ij}$ and the components of $s'_{\ell}$ with $\ell \le i+K$ or $\ell > j$ are kept fixed, 
we conclude that $s' \in B^{ij}$ is unique.

%Hence, once this subproblem is solved, 
%%if it is part of a larger DMDGP subproblem, 
%%it can be considered as a rigid part of a larger DMDGP problem containing it, 
%in order to keep the distance constraint associated to $\{i,j\}$ feasible, 
%positions $x_{i+K},\dots,x_j$ can only be reflected together, 
%meaning that we only need to care about its first symmetry hyperplane 
%when solving another subproblem that depends on $x_{i+K},\dots,x_j$. 
%%in order to compose with other partial reflections of the larger problem. 
%This will be explained in details in the next section. 

%\begin{remark}
%In view of Corollary~\ref{cor:1n}, if $x_i,\dots,x_j$ satisfies all distance constraints for $\{u,w\} \in E_D \cap E(G^{ij})$ 
%and $\|x_i - x_j \| = d_{ij}$, then $x_i,\dots,x_j$ satisfies all distance constraints for $\{u,w\} \in E(G^{ij})$, 
%i.e~it is a solution for $(G^{ij},K)$.
%\end{remark}

%%Notice that Theorem~\ref{thm:twos} implies that a fundamental $^K$DMDGP subproblem has only two solutions 
%%which are reflections of each other through the hyperplane defined by its first $K$ vertices. 
%%Hence, once this subproblem is solved and if it is part of a larger $^K$DMDGP subproblem, 
%%it can be handled as a {\it rigid body}, meaning that we only need to store its first symmetry hyperplane 
%%in order to compose with other partial reflections of the larger subproblem.

\section{New algorithm}\label{sec:algo}
Henceforth, we assume that subproblems $(G^{ij},K)$ spanned by pruning edges $\{i,j\}$ are solved following a given order $<$ in $E_P$ 
and that a realization $x \in \hat{X}$ is given. 
%and will allow us to generalize the concept of symmetry vertex set $S$ (see Eq.~\ref{eq:S}) to subproblems 
%$(G^{ij},K)$ spanned by pruning edges $\{i,j\}$. Moreover, 

\subsection{The conceptual algorithm}\label{sec:talgo}
First, we present a conceptual algorithm (Algorithm~\ref{algo:SBBUt}) 
which summarizes the ideas discussed in the previous sections. 
%We shall discuss a practical version in the next subsection. 

\begin{algorithm}[t]
\footnotesize
\begin{algorithmic}[1]
\STATE $\mathsf{SBBU}(G,K,(e_1,\dots,e_m), x \in \hat{X})$ 
\STATE Set $s=0$, $x(0) = x$
\FOR{$k = 1,2,\dots, m$}
	\STATE $\{i,j\} = e_k$
	\IF{$|S^{ij}|>0$}
		\STATE Find $s' \in B^{ij} \ : \ \| x_i(s') - x_j(s') \| = d_{ij}$ \label{st:finds} %\COMMENT{find position $x_j$} 
		\STATE Update $s=s'$ and $x(s) = U(x,s)$ \label{st:update}
	\ENDIF
\ENDFOR
\RETURN a valid realization $x(s)$
\end{algorithmic}
\caption{SBBU}  \label{algo:SBBUt}
\end{algorithm} 
%The other existing decomposition-based BP variants \cite{gramacho,fidalgo} are based on the given vertex order. The present proposal is based instead on an order on the pruning edges which is exploited in order to construct a realization recursively on the pruning edge spanned subproblems. As we shall see, this approach allows us much more solid theoretical and algorithmic treatments. As mentioned in Sect.~\ref{2}, each pruning edge $\{i,j\}$ induces a constraint on the set of valid partial reflection operations that can be applied to realizations of the vertices $i,\ldots,j$ (which we call the \textit{range} of the edge $\{i,j\}$). 

%in such a way that already solved DMDGP subproblems are considered as rigid parts of larger subproblems containing them.
%in which fundamental DMDGP subproblems are considered before larger subproblems containing them.
%Recall that each pruning edge $\{i,j\}$, with $i<j-K$, 
%defines a $^K$DMDGP subproblem $(G^{ij},K)$ that may be either   
%a fundamental $^K$DMDGP subproblem or a DMDGP subproblem containing many fundamental subproblems. 
%In the following subsections we describe how to solve and merge these subproblems to find a DMDGP solution. 

When solving subproblem $(G^{ij},K)$, if $|S^{ij}| = 0$ 
then this subproblem has already been solved implicitly, according to the following proposition.%(see Proposition~\ref{prop:SijZero}). 
\begin{proposition}\label{prop:SijZero}
Let $x(s)$ be a valid realization for $(G^{uw},K)$, for all $\{u,w\} \in P^{ij}$. 
If $S^{ij}  = \varnothing$, then $x(s)$ is valid for $(G^{ij},K)$.
\end{proposition}
\begin{proof}
If $ S^{ij}  = \varnothing$, then for every $v_{\ell} \in \{v_{i+K+1},\dots,v_j \}$, $\exists \{u,w\} \in P^{ij}$ 
such that $u+K+1 \le \ell \le w$. %$\ell \in [u+K+1,w]$. 
Suppose that $x(s)$ is such that $\| x_i(s) - x_j(s) \| \ne d_{ij}$. 
(a) By Theorem~\ref{thm:twos} and Assumption~\ref{a:feas}, $\exists s' \in \{0,1\}^{n-K}$ 
with some $s'_{\ell} \ne s_{\ell}$, for $\ell \in [u+K+1,w]$, such that $\| x_i(s') - x_j(s') \| = d_{ij}$. 
(b) Since $\ell \ge u+K+1$, it follows that $x_{w}(s') \notin \{ x_w(s), R_x^{u+K} (x_w(s)) \}$.
Thus, by Theorem~\ref{thm:twos}, $\| x_w(s') - x_u(s') \| \ne d_{uw}$. 
But (a) and (b) together contradict Assumption~\ref{a:feas}.
Hence $\| x_i(s) - x_j(s) \| = d_{ij}$ and the assertion follows from Lemma~\ref{lem:NaoEstraga}. 
\end{proof}

Otherwise, for $|S^{ij}|>0$, in Step~\ref{st:finds} we perform an exhaustive search to find $s' \in B^{ij}$ 
such that $\| x_i(s') - x_j(s') \| = d_{ij}$. In Step~\ref{st:update}, we update the current realization to $x(s')$ 
according to Eq.~\eqref{eq:xs}. 

\begin{theorem}\label{thm:correct}
Let $(G,K)$ be a feasible $^K$DMDGP instance. 
Considering exact arithmetic, Algorithm~\ref{algo:SBBUt} finds $x \in X$.
\end{theorem}
\begin{proof}
Since $x(0) = x \in \hat{X}$, 
due to Assumption~\ref{a:feas} and Lemma~\ref{lem:NaoEstraga}, Step~\ref{st:finds} is well-defined. 
From Remark~\ref{rem:key} and Step~\ref{st:finds}, it follows that $x(s') \in X(G^{ij})$, for every $e_k = \{i,j\}$. 
Thus, since for the last pruning edge $e_m = \{i_m,j_m\}$,  we have $E^{i_m,j_m} = E$, 
i.e~$G^{i_m,j_m}=G$, after this last subproblem is solved, $x(s) \in X(G^{i_m,j_m}) = X(G) = X$. 
\end{proof}

%%%%%%%%% Implementation %%%%%%%%%
\subsection{A practical algorithm}\label{sec:palgo}
In this section, based on a particular pruning edge order, 
we introduce a practical version of Algorithm~\ref{algo:SBBUt} which: 
\begin{enumerate}[i)]
\item does not required an initial realization $x \in \hat{X}$;
\item avoids the computation and storage of unnecessary reflectors $R_x^i(\cdot)$; 
\item may result in less operations in the update step (Step~\ref{st:update}) of Algorithm~\ref{algo:SBBUt};
\item allows us to discuss a concrete implementation for the sets $S^{ij}$.
\end{enumerate}

For this, instead of working with a full realization $x(s) \in \hat{X}$, 
which is updated through the binary vector $s$ by Eq.~\eqref{eq:xs}, 
and computing and storing reflectors $R_x^i(\cdot)$ based on $x=x(0) \in \hat{X}$, 
the idea is to grow a partial realization $x_1,\dots,x_t$, where $t = \arg \max \{ w \ | \ \{u,w\} \in P^{ij} \}$,  
and compute the necessary reflectors on the fly based on the current partial realization and $S^{ij}$. 
This way, for each subproblem $(G^{ij},K)$, we do not compute full valid realizations $x_1,\dots,x_n$ 
but valid partial realizations $x_1,\dots,x_j$, with $j \le t$.

\begin{assumption}\label{a:a2}
Pruning edges $\{i,j\}$, with $i < j$, 
are sorted in increasing order of $j$, followed by a decreasing order of $i$. 
\end{assumption}

Under this order, we can re-write the set of pruning edges preceding $\{i,j\}$ as 
\begin{equation}\label{eq:Pij0}
P^{ij} = \{ \{ u, w \} \in E_P \ | \ u < w < j \ \vee (w=j \wedge w>u>i)  \}.
\end{equation}

\begin{definition}\label{def:vpr}
We say that $x_1,\dots,x_t$, with $j \le t$, is a valid partial realization for $(G^{ij},K)$, 
if $x_1,\dots,x_t$ satisfies all distance constraints associated to edges in $\{ \{u,w\} \in E(G^{ij}) \ | \ w \le j\}$. 
\end{definition} 

\begin{remark}\label{rem:extend}
Recall that $E(G^{ij}) =  E^{ij} = E_D \cup P^{ij} \cup \{i,j\}$ and thanks to Assumption~\ref{a:a2}, 
there is no $\{u,w\} \in P^{ij}$, with $w>j$. This allows us to extend 
a valid partial realization $x_1,\dots,x_t$ for $(G^{ij},K)$ to a valid full realization $x_1,\dots,x_n$ for $(G^{ij},K)$, i.e~$x \in X(G^{ij})$, 
by simply growing $x_1,\dots,x_t$ to $x_1,\dots,x_n$ using discretization distances (see Subsection~\ref{sec:initx}),   
because no distance constraint $\{u,w\} \in P^{ij} \cup \{i,j\}$ is affected by this operation. 
\end{remark}

Assumption~\ref{a:a2}, along with \eqref{eq:Pij0}, will be used in the results that follow. 
Using this concepts we will show in the next subsections that 
when dealing with subproblem $(G^{ij},K)$:  
\begin{enumerate}
\item given a partial realization $x_1,\dots,x_t$ satisfying discretization distances and 
distances corresponding to pruning edges in $P^{ij}$, 
it can be extended to $x_1,\dots,x_j$ keeping feasibility of such distance constraints and new discretization constraints; 
\item it is possible to apply partial reflections to this extended partial realization in order to 
fulfill $\| x_i - x_j \| = d_{ij}$ without violating the distance constraints considered so far.
\end{enumerate}

\subsubsection{Initialization of candidate positions (Growth)}\label{sec:initx}
%In the previous sections we have discussed 
In Section~\ref{sec:solvesub} we shall explain how to find a valid partial realization for DMDGP subproblems $(G^{ij},K)$ 
by composing reflections through symmetry hyperplanes and applying them 
to positions $x_{i+K+1}, \dots, x_j$. %  of a current possible realization $x \in \hat{X}$. 
This procedure assumes that candidate positions for $x_i,\dots, x_j$ are available when we start to solve $(G^{ij},K)$. 
%It is important to use updated positions 
%%(which may change due to changes in previous positions $x_i, \dots, x_{i+K-1}$), 
%because the reflectors $R_x^{\ell}(\cdot)$ through symmetry hyperplanes are computed based on them. 
In this section, we describe how to initialize these positions. 

%%Before solving any pruning edge, we consider $x$ as any possible realization in the search tree: 
%%recall from Definition~\ref{def:dmdgp} that the first $K$ vertices can be localized uniquely (up to rotations and translations), 
%%and if we consider only discretization distances, the other vertices $v_i$, with $i>K$ in the discretization order, 
%%can be localized  by solving the intersection of $K$ spheres in $\R^K$ which gives two possible positions $\{x_i^{-}, x_i^{+}\}$. 
%%Let $x$ be the initial realization obtained by choosing $x_i^{-}$, for every $i=K+1,\dots, n$, for example. 
%%
%%Now, every time we start to solve $(G^{ij},K)$, there is a current realization $x$, 
%%that can be the initial realization we just mentioned, or an updated version of it. 

%Of course, in principle, two subproblems $(G^{ij},K)$ and $(G^{uv},K)$ that share at least $K+1$ vertices  
%could be solved separately, each one in its own coordinate system, and only after merging 
%(e.g., by a Procrustes analysis on their common $K+1$ vertices) to find the absolute position of their vertices:  
%a process known by {\it synchronization}. Two or more subproblems that share less than $K+1$ vertices 
%can have the global position of their vertices decided only in view of pruning edges covering them (which shall be resolved later). 

%In order to avoid this synchronization process, and directly compute all vertex positions in a global coordinate system, 

From now on, we assume that initialization of candidate positions must follow the vertex order from $v_1$ to $v_n$, 
meaning that if $(G^{ij},K)$ is the current subproblem, and $x_t$ is the last initialized position, such that $t < j$, 
then we initialize positions from $x_{t+1}$ to $x_j$, whereas $x_1,\dots, x_t$ remain unchanged. 
In other words, the candidate positions $x_{t+1},\dots,x_j$ 
are grown from the current partial realization $x_1,\dots,x_t$ using only distance constraints associated to discretization edges.  
Moreover, each position $x_i$ is initialized only once, although it can be modified later (see Section~\ref{sec:solvesub}) 
in order to satisfy a distance constraint corresponding to a pruning edge $e' \ge e = \{i,j\}$. 
This is formalized in Proposition~\ref{prop:initx}.

\begin{proposition}\label{prop:initx}
%If DMDGP subproblems are solved following the above pruning edge order, 
Assume edges in $E_P$ are ordered as $e_1,\dots,e_m$. 
Then, before solving $(G^{ij},K)$, positions $x_1,\dots, x_j$  can be initialized such that
\begin{equation}\label{eq:ED}
\forall \{ \ell, k \} \in E_D \cap E(G[v_1,\dots,v_j]), \quad \| x_{\ell} - x_k \| = d_{\ell k},
\end{equation}
and
\begin{equation}\label{eq:EPij}
\forall \{ \ell, k \} \in P^{ij}, \quad \| x_{\ell} - x_k \| = d_{\ell k}.
\end{equation}
\end{proposition}
\begin{proof}
We prove this by induction on the edge order. 
In the base case we consider $e_1 = \{i_1,j_1\} \in E_P$ spanning the first subproblem to be solved. 
The positions $x_i$, for $i=1,\dots,j_1$ are initialized right away.  
From Definition~\ref{def:dmdgp}, $x_1,\dots,x_K$ can be localized uniquely (up to rotations and translations) 
by different methods \cite{dokmanic2015,alencar2019}. Hence, $\| x_{\ell} - x_k \| = d_{\ell k}, \forall \{\ell, k\} \in E(G[v_1,\dots,v_K])$. 
%Then, by considering only discretization distances, each vertex $v_i$, with $K<i \le j_1$ in the discretization order, 
%can be localized  by finding the intersection of $K$ spheres in $\R^K$ which gives two possible positions $\{x_i^{-}, x_i^{+}\}$. 
Then, by $K$-lateration (see Remark.~\ref{rem:Klat}), there are at most two positions $\{x_i^+,x_i^-\}$ for $v_i\in V$ for each $K<i\le j_1$.
Notice that any partial realization $x_1,\dots,x_{j_1}$ is enough to build partial reflections. 
The correct alternative will be chosen later by the appropriate partial reflection composition which 
satisfies constraints defined by pruning edges (Section~\ref{sec:solvesub} gives more details). 
Thus, without loss of generality, let $x_1, \dots, x_{j_1}$ be the partial realization obtained by choosing $x_i^{-}$, for every $i=K+1,\dots, j_1$. Since $P^{i_1 j_1} = \varnothing$, this partial realization satisfies \eqref{eq:ED} and \eqref{eq:EPij} for $\{i,j\} = \{i_1,j_1\}$.

The induction hypothesis is that \eqref{eq:ED} and \eqref{eq:EPij} hold for pruning edges $e_1,\dots,e_k$, 
i.e~$x_1,\dots,x_t$ is a valid partial realization for all subproblems spanned by these edges, where 
$$
t:= \max \{ w \ | \ \{u,w\} \in P^{ij} \} = \max \{ w \ | \ \{u,w\} \in \{e_1,\dots,e_k\} \}.
$$ 

In the inductive step, let us prove that \eqref{eq:ED} and \eqref{eq:EPij} also hold for pruning edge $e_{k+1} = \{i,j\}$ 
spanning subproblem $(G^{ij},K)$. 
%Let $(G^{ij},K)$ be the next subproblem to be solved, spanned by pruning edge $e_{k+1} = \{i,j\}$. 

Since subproblems spanned by edges in $P^{ij}$ are solved, positions $x_1,\dots,x_t$ are already initialized 
and satisfy \eqref{eq:EPij}, and \eqref{eq:ED} with $v_j = v_t$. 

If $j \le t$, then there is nothing left to do. 
Thus, suppose $j>t$. Then, positions $x_{t+1},\dots,x_{j}$ can be initialized 
%by a sequence of sphere intersections (see Algorithm~\ref{algo:initx}) 
by $K$-lateration (see Remark~\ref{rem:Klat} and Algorithm~\ref{algo:initx}) 
based on discretization distances such that \eqref{eq:ED} holds. 
\end{proof}
%and we proceed as will be described in Section~\ref{sec:solvemerged} to solve $(G^{i_{\ell+1}j_{\ell+1}},K)$. 
\begin{remark}\label{rem:nosync}
The proof of Proposition~\ref{prop:initx} describes a procedure for initialization of $x_1,\dots,x_j$ before solving $(G^{ij},K)$. 
It is important to notice that such initialization is done sequentially and depends on previously 
computed positions which are not recomputed in this step. Thus, after the initialization, the current partial realization 
continues to be valid for all already solved subproblems. 
%Moreover, due to this sequential way of computing positions, no post-synchronization is necessary to stitch realizations of subproblems together. 
\end{remark}

\begin{algorithm}[t]
\footnotesize
\begin{algorithmic}[1]
\STATE $\mathsf{InitializePositions}(x,t,j)$
\IF{$t \ge j$}
	\RETURN $x$ and $t$.
\ENDIF
\IF{$t=0$}
	\STATE Initialize $x_1,\dots, x_K$ as a solution of $\| x_{\ell} - x_i \|^2 = d_{\ell i}^2, \forall \{\ell,i\} \in E(G[v_1,\dots,v_K])$
	\STATE Set $t=K$
\ENDIF
\FOR{$i = t+1,\dots, j$}
	      \STATE Find solutions $\{x_i^+, x_i^{-}\}$ for the system: $\| x_{\ell} - x_i \|^2 = d_{\ell,i}^2, \ell = i-K,\dots,i-1.$ %Compute the intersection of spheres centered at $x_{i-K},\dots,x_{i-1}$, with radii $d_{i-K,i},\dots,d_{i-1,i}$ 
%      to obtain the candidate positions $\{x_i^+, x_i^{-}\}$. 
      \STATE Set $x_i = x_i^{-}$
\ENDFOR
\STATE Set $t=j$
\RETURN $x$ and $t$.
\end{algorithmic}
\caption{InitializePositions}  \label{algo:initx}
\end{algorithm} 

Algorithm~\ref{algo:initx} gives a pseudocode for the function $\mathsf{IntializePositions}$
which receives a current realization $x$ (actually, the current partial realization $x_1,\dots,x_t$), 
the index of the last initialized position $t$, 
and the index $j$ of the last vertex whose position needs initialization. 
Updated $x$ and $t$ are returned by this function.

\subsubsection{Solving DMDGP subproblems (Correction)}\label{sec:solvesub}
Now we explain how to find a valid partial realization for a DMDGP subproblem $(G^{ij},K)$, 
given a valid partial realization $x_1,\dots,x_t$  %with $t \ge j$, 
for $(G^{uw},K), \ \forall\{u,w\} \in P^{ij}$. 
%a striking observation is that we do not need to recompute {\it all} the positions $x_{i+K+1},\dots,x_j$ 
%every time $d_{ij}$ is violated. 
%every time we end up in an infeasible node in the corresponding search tree. 

Recall from Proposition~\ref{prop:reflect} that given positions $x_{i+1},\dots,x_{i+K+1},\dots,x_j$, {\it valid or not},  
%(e.g.~we can always take the leftmost path in the corresponding search tree),   
we can build all necessary symmetry hyperplanes and their corresponding reflection operators $R_x^{i+K+1}(\cdot), \dots, R_x^{j}(\cdot)$. 
Then, based on Theorem~\ref{thm:twos} and Proposition~\ref{prop:two},  
we can apply compositions of such reflection operators {\it only} to $x_j$ 
until we find its correct position, as illustrated in Figure~\ref{fig:reflections} for the 2D case. 

The only decision to be taken is whether each of the reflectors $R_x^{i+K+1}(\cdot)$, $\dots$, $R_x^{j}(\cdot)$ should be applied or not to $x_j$ 
in order to fulfill $\| x_i  - R(x_j,\bar{s}) \| = d_{ij}$, where $R(.,\bar{s})$ is a composition of the chosen reflectors:
\begin{equation}\label{eq:compor}
R(y,\bar{s}) := (R_x^{i+K+1})^{\bar{s}_1} (R_x^{i+K+2})^{\bar{s}_2}\dots (R_x^{j})^{\bar{s}_{j-i-K}}(y).
\end{equation}
In \eqref{eq:compor}, the binary vector $\bar{s}$ is of size $j-i-K$, and $(R_x^{\ell})^0 := I$, the identity operator in $\R^K$, 
whereas $(R_x^{\ell})^1 := R_x^{\ell}$, for $\ell = i+K+1, \dots, j$. 

In contrast to Algorithm~\ref{algo:SBBUt}, where $s$ is the global binary decision variable and all the reflectors 
are computed based on the first realization $x(0) = x \in \hat{X}$, 
now the reflectors $R_x^{i+K+\ell}(\cdot)$ for which $v_{i+K+\ell} \in S^{ij}$ 
are computed based on the current partial realization $x_1,\dots,x_t$, 
and $\bar{s}$ is a local binary decision variable belonging to
\begin{equation}\label{eq:barBij}
\bar{B}^{ij} := \{ \bar{s} \in \{0,1\}^{j-i-K} \ | \ \bar{s}_{\ell} = 0 \text{ if } v_{i+K+\ell} \notin S^{ij} \}. 
\end{equation}

%More specifically, to compute a realization for a subproblem $(G^{ij},K)$ given positions initialized as in Proposition~\ref{prop:initx}, 
Thus, we look for a binary vector $\bar{s} \in \bar{B}^{ij}$ 
%We look for a \tred{binary vector} 
such that 
\begin{equation}\label{eq:compor2}
x'_j = R(x_j,\bar{s}) = (R_x^{i+K+1})^{\bar{s}_1} (R_x^{i+K+2})^{\bar{s}_2}\dots (R_x^{j})^{\bar{s}_{j-i-K}}(x_j) 
\end{equation}
satisfies $\| x_i - x'_j\| = d_{ij}$.  We remark that this search is exhaustive: we test all $|\bar{B}^{ij}| = 2^{|S^{ij}|}$ possible choices for $\bar{s}$ 
(recall that there is a unique $\bar{s}$ that works, as discussed after Proposition~\ref{prop:two}).

Once $\bar{s}$ is found, the positions of $v_{i+K+1},\dots,v_{t}$ are updated according to:
\begin{equation}\label{eq:reflectall2}
%x'_{\ell} = \bar{U}(x_{\ell},\bar{s}) := \begin{cases}
%x_{\ell}, & \quad \ell = i,i+1,\dots,i+K \\
%\left( \prod_{t=1}^{\ell-i-K} (R_x^{i+K+t})^{\bar{s}_t} \right) (x_{\ell}),	& \quad \ell \ge i+K+1. %\ell = i+K+1,\dots,j.
%\end{cases}
x'_{\ell} = \bar{U}(x_{\ell},\bar{s}) : = \left( \prod_{t=1}^{\ell-i-K} (R_x^{i+K+t})^{\bar{s}_t} \right) (x_{\ell}), \quad  \ell = i+K+1,\dots,t.
\end{equation}
This update maintain feasibility of $x_1,\dots,x_t$ with respect to $(G^{uw},K)$, for every $\{u,w\} \in P^{ij}$, 
because positions $x_{u+K},\dots,x_w$ 
%can only be changed all-together by 
%must all be updated by the partial reflections 
are only updated simultaneously by partial reflections 
$R_{x}^{v}(x_{u+K}), \dots, R_{x}^{v}(x_w)$, for $v \le u+K$.

\subsubsection{Symmetry vertex sets}
\begin{algorithm}[t]
\footnotesize
\begin{algorithmic}[1]
\STATE $\mathsf{SBBU}(G,K)$ 
\STATE Order edges $\{v_i,v_j\} \in E_P$ in increasing order of $j$ and decreasing order of $i$ obtaining a sequence $(e_1,\dots,e_m)$, with $m = |E_P|$. Set $t = 0$, $n=|V|$ 
\STATE Set ${\cal C} = \{\{v_i\}\}_{i=K+1}^n$ 
\STATE $\mathsf{InitializePositions}(x,t,K)$ \COMMENT{positions for the initial clique}
\FOR{$k = 1,2,\dots, m$}
	\STATE $\{i,j\} = e_k$
	\IF{$\rho_{\cal C}(i+K) \ne \rho_{\cal C}(j)$}
		\STATE $\mathsf{InitializePositions}(x,t,j)$ \label{s2:initxi}
		\STATE Set $C^0 = \rho_{\cal C}(i+K)$ and ${\cal D} = \rho_{\cal C}(\{ i+K+1, \dots, j \})\setminus \{ C^0 \}$  \label{s2:subC}
%%		\STATE Set ${\cal D} = \{ \rho_{\cal C}( \ell) \ : \ i + K < \ell \le j \}$ and $C^0 = \rho_{\cal C}(i+K)$ \label{s2:subC}
%%		\STATE Set ${\cal D} = \cup_{i + K < \ell \le j} \ \rho_{\cal C}(\ell)$ and $C^0 = \rho_{\cal C}(i+K)$
		\STATE Let $S^{ij} = \cup_{C \in {\cal D}} \mathsf{first}(C)$ \COMMENT{local symmetry vertices}  \label{s2:planes}
		\STATE Compute $R_x^{\ell}(\cdot)$ for each $v_{\ell} \in S^{ij}$\label{s2:planes_end}
		\STATE Find $s \in \bar{B}^{ij} \ : \ \| x_i - R(x_j,\bar{s}) \| = d_{ij}$ \COMMENT{find position $x_j$} \label{s2:finds}
%%		\STATE If such $s$ does not exist, \textbf{return} ``infeasible instance''
		\STATE Update $x_{\ell} = \bar{U}(x_{\ell},\bar{s})$, for $\ell = i+K+1,\dots,j$\label{s2:reflectall}
		\STATE Set $C^{+} = \left( \cup_{C \in {\cal D}} C \right) \cup C^0$ and update ${\cal C} = ({\cal C} \setminus ( {\cal D} \cup \{C^0\}) ) \cup C^{+}$ \label{s2:merge}%\COMMENT{merge} 
	\ENDIF
\ENDFOR
\IF{$t<n$}
	\STATE $\mathsf{InitializePositions}(x,t,n)$ \label{s2:remaining}
\ENDIF
%%\STATE Set $S = \cup_{C \in {\cal C}} \mathsf{first}(C)$ \COMMENT{symmetry vertices} \label{s2:setS}
\RETURN a valid realization $x$ %and the set of symmetry vertices $S$
\end{algorithmic}
\caption{SBBU}  \label{algo:SBBU2}
\end{algorithm} 
The ideas of the Subsections~\ref{sec:initx} and \ref{sec:solvesub} lead to Algorithm~\ref{algo:SBBU2}. 
This algorithm makes use of ${\cal C}$, a partition of $\{v_{K+1}, \dots, v_n\}$ 
used to obtain the sets $S^{ij}$. 
At the beginning, we set ${\cal C} = \{ \{ v_i \} \}_{i=K+1}^n$. 
This partition is updated in Step~\ref{s2:merge} taking into account 
%subproblems that are dependent overlapping with $(G^{ij},K)$ and the vertices they cover. 
already solved subproblems. 
Assume that subsets of vertices in ${\cal C}$ are ordered according to the vertex order of Definition~\ref{def:dmdgp}. 
Let us denote by $\mathsf{first}(C^0)$ the first vertex of $C^0 \in {\cal C}$. 
%We remark that the sets in ${\cal C}$ are always disjoint and their union is $\{v_{i+K},\dots,v_{n}\}$. 

We also introduce a function $\rho_{\cal C}: \{K+1,\dots,n\}  \rightarrow {\cal C}$, 
parametrized by ${\cal C}$, such that $\rho_{\cal C}(\ell)$ returns the unique element of ${\cal C}$ containing vertex $v_{\ell}$. 
The next proposition shows that this function is well-defined at every iteration of Algorithm~\ref{algo:SBBU2}. 

\begin{proposition}\label{prop:partition}
At every iteration of Algorithm~\ref{algo:SBBU2}, ${\cal C}$ is a partition of the subset of vertices $\{ v_{K+1}, \dots, v_n \}$. 
\end{proposition}
\begin{proof}
At the first iteration ${\cal C} = \{ \{ v_i \} \}_{i=K+1}^n$. 
Assume ${\cal C}$ is a partition of $\{ v_{K+1}, \dots, v_n \}$ at the beginning of iteration $k$. 
If $\rho_{\cal C}(i+K) = \rho_{\cal C}(j)$, then we go to the next iteration with ${\cal C}$ unchanged. 
Otherwise, from Step~\ref{s2:subC}, ${\cal D}$ and $C^0$ are subsets of ${\cal C}$. 
Then, Step~\ref{s2:merge} updates ${\cal C}$ by removing these subsets and including their union, 
hence, the updated ${\cal C}$ is still a partition of $\{ v_{K+1}, \dots, v_n \}$. 
\end{proof}

\begin{proposition}\label{prop:C}
In Algorithm~\ref{algo:SBBU2}, after Step~\ref{s2:merge}, there exists a unique $C \in {\cal C}$ 
such that $C \supset \{v_{i+K},\dots,v_{j} \}$.
\end{proposition}
\begin{proof}
Follows directly from Steps~\ref{s2:subC} and \ref{s2:merge} of Algorithm~\ref{algo:SBBU2}. 
\end{proof}

%Recall that $\rho_{\cal C}(\ell)$ returns the unique element of ${\cal C}$ containing $v_{\ell}$. 
We remark that $\rho_{\cal C}(\{ i+K+1, \dots, j \})$ denotes the image of $\{ i+K+1, \dots, j \}$ by $\rho_{\cal C}$ 
in the definition of ${\cal D}$ (Step~\ref{s2:subC}), 
i.e~it returns elements of ${\cal C}$ 
whose union contains $v_{i+K+1},\dots,v_j$ and $C^0 = \rho_{\cal C}(i+K)$ is the element of ${\cal C}$ containing $v_{i+K}$.

\begin{proposition}\label{prop:SijEmpty}
In Algorithm~\ref{algo:SBBU2}, $\rho_{\cal C} (i+K) = \rho_{\cal C}(j)$ if and only if $S^{ij} = \varnothing$.
\end{proposition}
\begin{proof}
If $S^{ij} = \varnothing$, then 
\begin{equation}\label{eq:key}
\forall \ v_{\ell} \in \{ v_{i+K+1},\dots,v_j \}, \quad \exists\{u,w\} \in P^{ij} \ : \ u+K < \ell \le w \le j.
\end{equation}
In particular, for $\ell = i+K+1$, there exists $\{ r, z \} \in P^{ij}$ such that
$r + K < i+K+1 \le z \le j$. Clearly, $r \le i$. 
From Proposition~\ref{prop:C} there exists a unique $C^1 \in {\cal C}$ 
such that $C^{1} \supset \{ v_{r+K},\dots, v_{i+K},v_{i+K+1},\dots,v_z\}$. 

Thus, if $z=j$, then $\rho_{\cal C}(i+K) = \rho_{\cal C}(j)$. 

Otherwise, for $z < j$, because $i+K+1 \le z$, it follows that $v_{z+1} \in \{ v_{i+K+1},\dots,v_j \}$ and from \eqref{eq:key}, 
there exists $\{u,w\} \in P^{ij}$ such that $u+K < z+1 \le w \le j$ (clearly, $u + K \le z$).
Thus, from Proposition~\ref{prop:C}: 
\begin{equation}\label{eq:c2}
\exists ! C^2 \supset \{ v_{u+K}, \dots, v_z, v_{z+1},\dots,v_w \}.
\end{equation}

If $u \le r \le i$, then from \eqref{eq:c2}, we obtain $v_{i+K} \in C^2$. 

Otherwise, for $r < u \le z - K$, then $r+K < u+K \le z$, implying that $v_{u+K} \in C^1$.
In either case, we have $\rho_{\cal C}(u+K) = \rho_{\cal C} (i+K)$. 
From Proposition~\ref{prop:partition},  we conclude that  $\rho_{\cal C}(i+K) = \rho_{\cal C}(w)$. 

Hence, if $w=j$, $\rho_{\cal C}(i+K) = \rho_{\cal C}(j)$. 

Otherwise ($w<j$), in view of \eqref{eq:key}, we can apply the same argument to $v_{w+1}$, 
and repeat until we find $\{h,p\} \in P^{ij}$ with $p=j$.

On the other hand, to prove that $\rho_{\cal C}(i+K) = \rho_{\cal C}(j)$ implies $S^{ij} = \varnothing$, we use the counter-positive.
Suppose there exists $v_{\ell} \in \{v_{i+K+1},\dots,v_j \}$ such that $\not\exists \{u,w\} \in P^{ij}$ such that $u+K < \ell \le w \le j$. 
This means that $\forall \{u,w\} \in P^{ij}$ either (i) $w < \ell$ or (ii) $\ell \le w \le j$ and $u+K \ge \ell$. 
If $w < \ell, \forall \{u,w\} \in P^{ij}$, then $\rho_{\cal C}(j) = \{v_j\} \ne \rho_{\cal C}(i+K)$, because $w < \ell \le j$ ($v_{\ell}$ and $v_j$ were never reached). 
 
Thus, let us consider $\{u,w\} \in P^{ij}$ such that $\ell \le w \le j$ and $u+K \ge \ell$. 
Without loss of generality, assume $w=j$. 
Since $\ell \ge i+K+1$, then $u+K \ge i+K+1$ (or $u \ge i+1$), implying that $\rho_{\cal C}(i+K) \ne \rho_{\cal C}(u+K) = \rho_{\cal C}(j)$, 
where the last equality follows from Proposition~\ref{prop:C}. 
\end{proof}

%When solving subproblem $(G^{ij},K)$, 
Proposition~\ref{prop:SijEmpty} shows that if $v_{i+K}$ and $v_j$ are in the same subset, i.e., $\rho_{\cal C}(i+K) = \rho_{\cal C}(j)$, 
then this subproblem was already solved implicitly (see Proposition~\ref{prop:SijZero}). % \tred{due to the merging of dependent overlapping subproblems:}
This is equivalent to condition $|S^{ij}| = 0$ in Algorithm~\ref{algo:SBBUt}. 

%In Step~\ref{s2:initxi}, we initialize possible positions for $x_{t+1},\dots, x_j$, 
%where $t$ stores the index of the last initialized position (see Section~\ref{sec:initx}). 
Otherwise, we need to obtain the set $S^{ij}$ of symmetry vertices for $(G^{ij},K)$. 
This is accomplished in Step~\ref{s2:planes}. % where $\mathsf{first}(C)$ retrieves the first vertex of subset $C$ (w.r.t. vertex order of Definition~\ref{def:dmdgp}).

\begin{theorem}\label{thm:SijOK}
If $\rho_{\cal C}(i+K) \ne \rho_{\cal C}(j)$, then $\bigcup_{C \in {\cal D}} \mathsf{first}(C) = S^{ij}$.
\end{theorem}
\begin{proof}
Let $v_{\ell} \in \cup_{C \in {\cal D}} \first{C}$. 
From Proposition~\ref{prop:partition}, $\exists! \hat{C} \supset \{ v_{\ell} \}$ such that $\hat{C} \in {\cal D}$ and $\first{\hat{C}} = v_{\ell}$. 
Suppose $v_{\ell} \not\in S^{ij}$, i.e~there exists $\{u,w\} \in P^{ij}$ such that $u+K < \ell \le w$. 
From Proposition~\ref{prop:C}, $\exists! C \in {\cal C}$ such that $C \supset \{v_{u+K},\dots,v_w \}$ and since ${\cal C}$ 
is a partition, it follows that $C = \hat{C}$. But $\first{\hat{C}} = \first{C} = v_{u+K} \ne v_{\ell}$ contradicting $\first{\hat{C}} = v_{\ell}$. 
Therefore, $\not\exists \{u,w\} \in P^{ij}$ such that $u+K < \ell \le w$. Thus, $v_{\ell} \in S^{ij}$. 

Conversely, let $v_{\ell} \in S^{ij}$. Then, for every $\{u,w\} \in P^{ij}$ either (i) $w < \ell $ or (ii) $\ell \le w \le j$ and $u+K \ge \ell$. 
If $\forall \{u,w\} \in P^{ij}$, we have $w < \ell$, then $C = \rho_{\cal C}(\ell) = \{v_{\ell}\} \in {\cal D}$ and $\first{C} = v_{\ell}$. 
Otherwise, there are $\{u,w\} \in P^{ij}$ such that $\ell \le w \le j$. For all of these, $u+K \ge \ell$. 
Recall from Algorithm~\ref{algo:SBBU2} that $\first{C^+} = \first{C^0} = \first{\rho_{\cal C}(u+K)}$. 
We split the analysis in three cases.

{\it Case 1}: $v_{\ell} \not\in \rho_{\cal C}(u+K)$. Then, $\ell < u+K$, implying that $v_{\ell} < \first{u+K}$.
%Thus, after solving $(G^{uw},K)$, we have $\rho_{\cal C}(\ell) = \{v_{\ell}\}$. 
Thus, after iteration $k$ with $e_k = \{u,w\}$, we have $\rho_{\cal C}(\ell) = \{v_{\ell}\}$. 

{\it Case 2}: $v_{\ell} \in \rho_{\cal C}(u+K)$ and $\first{\rho_{\cal C}(u+K)} = v_{\ell}$. In this case, 
%after solving $(G^{uw},K)$, 
after iteration $k$ with $e_k = \{u,w\}$, 
$\rho_{\cal C} (\ell) = \rho_{\cal C} (u+K) = \{ v_{\ell}, \dots, v_{u+K}, \dots, v_p \}$. 
Thus $\first{\rho_{\cal C}(\ell)} = v_{\ell}$.

{\it Case 3}: $v_{\ell} \in \rho_{\cal C}(u+K)$ but $\first{\rho_{\cal C}(u+K)} < v_{\ell}$. 
In this case, $C^0 = \rho_{\cal C}(\ell) = \rho_{\cal C}(u+K) = \{ v_q, \dots, v_{\ell},\dots,v_{u+K},\dots,v_p \}$. 
The set $C^0$ 
% is the result of solving $(G^{hp},K)$, 
is the result of iteration $k'$, with $e_{k'} = \{h,p\} < \{u,w\} = e_k$. Clearly $\{h,p\} \in P^{ij}$.  
%for $\{h,p\} \in P^{ij}$ such that $\{h,p\} < \{u,w\}$. 
From Proposition~\ref{prop:C}, $\exists ! C \supset \{ v_{h+K},\dots,v_p \}$. 
Notice that $u+K \le p \le w$. 
Since $q > h+K$ contradicts the fact that $\first{\rho_{\cal C}(u+K)} = v_q$, 
then $q \le h+K < \ell$. This leads to $h+K < \ell \le u+K \le p$ which implies that $v_{\ell} \not\in S^{ij}$, a contradiction. 

Therefore, only cases 1 and 2 can happen and both imply in $v_{\ell} \in \bigcup_{c \in {\cal D}} \mathsf{first}(C)$. 
\end{proof}

\begin{corollary}\label{thm:correctp}
Let $(G,K)$ be a feasible $^K$DMDGP instance. 
Considering exact arithmetic, Algorithm~\ref{algo:SBBU2} finds a valid realization $x$ for $(G,K)$. 
%Furthermore, the algorithm correctly returns the set $S$ of symmetry vertices of $(G,K)$. 
%or detects that $(G,K)$ is infeasible. 
\end{corollary}
\begin{proof}
%The first part follows from Lemma~\ref{lem:SijOK} and Theorem~\ref{thm:correct}. 
From Theorem~\ref{thm:SijOK} and Remark~\ref{rem:extend} correctness of Algorithm~\ref{algo:SBBU2} 
follows from Theorem~\ref{thm:correct}.
\end{proof}

In the end, we obtain a valid partial realization $x_1,\dots,x_t$ for all subproblems $(G^{ij},K)$, 
with $\{i,j\} \in E_P$. If $t=n$, we are done. Otherwise, in view of Remark~\ref{rem:extend}, 
$x_1,\dots,x_t$ can be extended to a valid realization $x \in X$. 
This explains Step~\ref{s2:remaining}. 

%Since each subproblem $(G^{ij},K)$ spanned by $\{i,j\} \in E_P$ is correctly solved (Lemma~\ref{lem:GijOK}) 
%and its solution does not change the feasibility of $x_1,\dots, x_j$ for previous subproblems (Lemma~\ref{lem:NaoEstraga}), 
%then after the last subproblem, spanned by $\{i_m,j_m\}$, is solved, we have that $x_1,\dots,x_t$, with $t=j_m$, 
%satisfies \linebreak$\| x_i - x_j \| = d_{ij}, \forall \{i,j\} \in E(G^{1t})$. If $t=n$, then we have a valid realization for \eqref{eq:dgp}. 
%Otherwise, if $t<n$, we initialize the remaining positions $x_{t+1},\dots,x_n$ respecting the remaining discretization distances. 
%This proves the first statement. 

%\tred{For the second part, we need to prove that 
%$$
%\cup_{C \in {\cal C}} \mathsf{first}(C) = S:= \{v_{\ell} \in V:\nexists \{v_i,v_j\}\in E\text{ with }i+K<\ell\leq j\}
%$$
%holds at Step~\ref{s2:setS}. This can be done by applying Lemma~\ref{lem:SijOK} to a hypothetical 
%subproblem $(G^{1n},K)$ spanned by pruning edge $\{1,n\}$.}

%%% floating-point arithmetic
Even under Assumption~\ref{a:feas}, due to floating point arithmetic, 
in Step~\ref{s2:finds} we may not be able to find $s$ such that $\| x_i - R(x_j,\bar{s})\| = d_{ij}$. 
Thus, instead of stopping as soon as we find a $s$ such that $| \| x_i - R(x_j,\bar{s})\| - d_{ij}| \le \varepsilon$, 
for a prescribed tolerance $\varepsilon$, 
we actually consider all $2^{|S^{ij}|}$ possibilities and 
choose $s$ for which $| \| x_i - R(x_j,\bar{s})\| - d_{ij}|$ is minimum. 
In case $| \| x_i - R(x_j,\bar{s})\| - d_{ij}| > \varepsilon$ for every $\bar{s} \in \bar{B}^{ij}$, 
then we actually interrupt the algorithm and return ``failure''. 
However, this never happened in the numerical experiments of Section~\ref{sec:results}. 

%%% computational cost

We remark that  $|S^{ij}|$ is a good indicator of the computational cost for solving subproblem $(G^{ij},K)$,  
because it determines the number of $2^{|S^{ij}|}$ reflection compositions that we need to apply to $x_j$ 
in order to find its correct position.  Thus, we define the corresponding {\it total work} 
to solve a $^K$DMDGP instance as
\begin{equation}\label{eq:work}
W  := \sum_{\{i,j\} \in \hat{E}} 2^{|S^{ij}|}, 
\end{equation}
where $\hat{E} = \{ \{v_i,v_j\} \in E_P \ | \ |S^{ij}| > 0 \}$. 
Let us also denote by $\bar{W} = \max_{\hat{E}} 2^{|S^{ij}|}$, 
the maximum work per pruning edge. 
We also remark that $|S^{ij}|$ depends on the order in which the pruning edges are handled 
and, in this paper, we consider only the order described in Assumption~\ref{a:a2}.

\section{Computational results}\label{sec:results}
An efficient implementation of the function $\rho_{\cal C}$ needs to deal with its evaluation and the update of the subsets of ${\cal C}$.
We adopted the structure proposed by Newman and Ziff \cite{newman2001}, 
which allows the evaluation of $\rho_{\cal C}$ and the subsets update in time $O(\log_2(|V|)) $ and memory $O(|V|)$.

In order to validate Algorithm~\ref{algo:SBBU2} and assess its performance,  
we generate a set of protein-like instances ($K=3$) whose data were extracted from Protein Data Bank (PDB) \cite{pdb}, 
and compare the results with those of the classic BP \cite{llm08,coap2012}. 
For each protein, we consider only the backbone composed by the sequence of atoms $N-C_{\alpha}-C$ 
and include an edge in the corresponding graph:
\begin{enumerate}
\item either when the atoms are separated by at most three covalent bonds 
\item or the distance between pairs of atoms is smaller than a certain cut-off value. 
\end{enumerate}
The resolution of NMR experiments usually varies between 5 \AA \ and 6 \AA. 
The smaller the cut-off value, the sparser the DMDGP instance.  
Each instance was generated by the first model and first chain of the PDB file.

\begin{figure}
\centering
\includegraphics[scale=0.3,trim=100 0 100 0,clip]{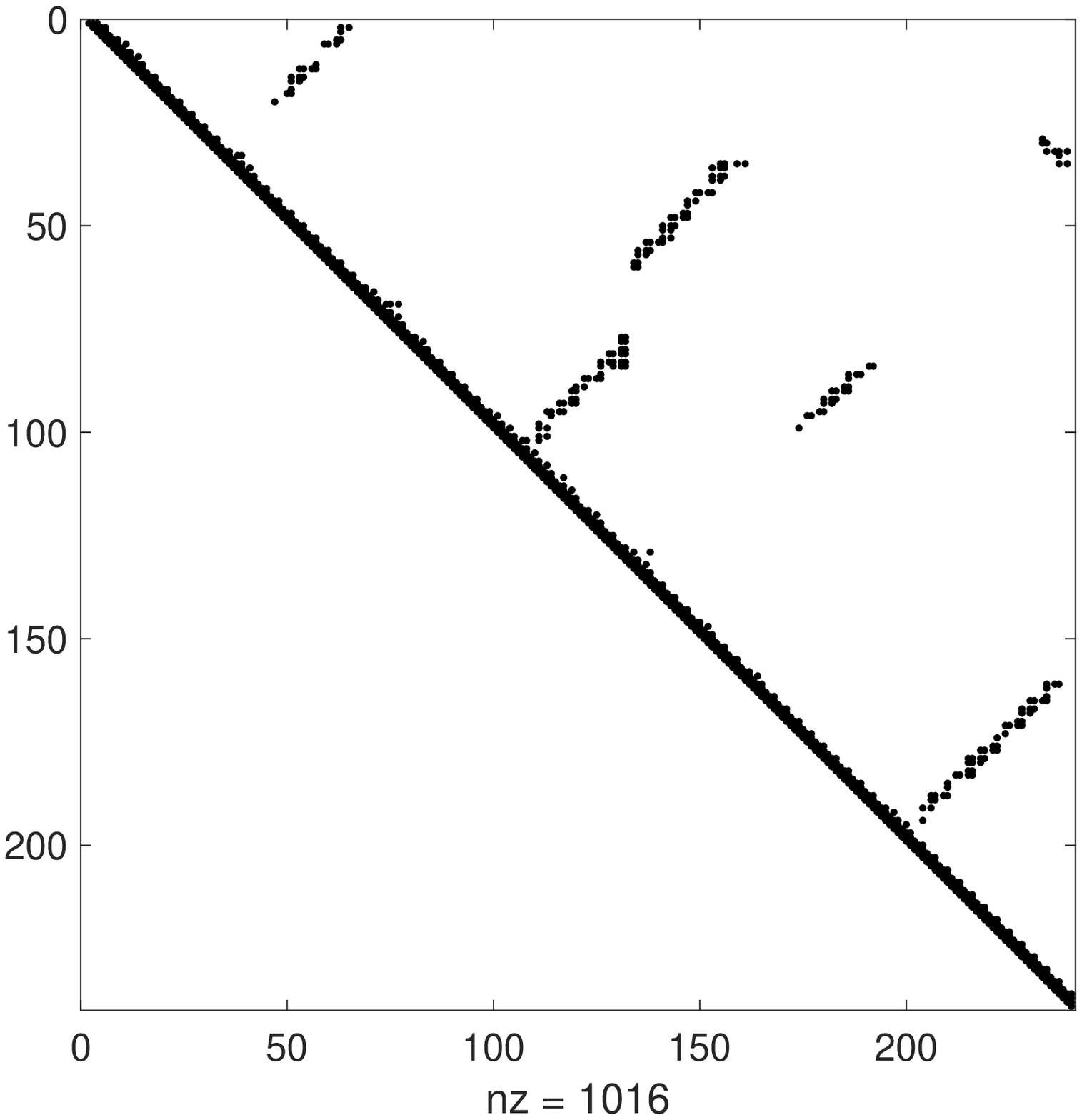}
\includegraphics[scale=0.305,trim=100 0 100 0,clip]{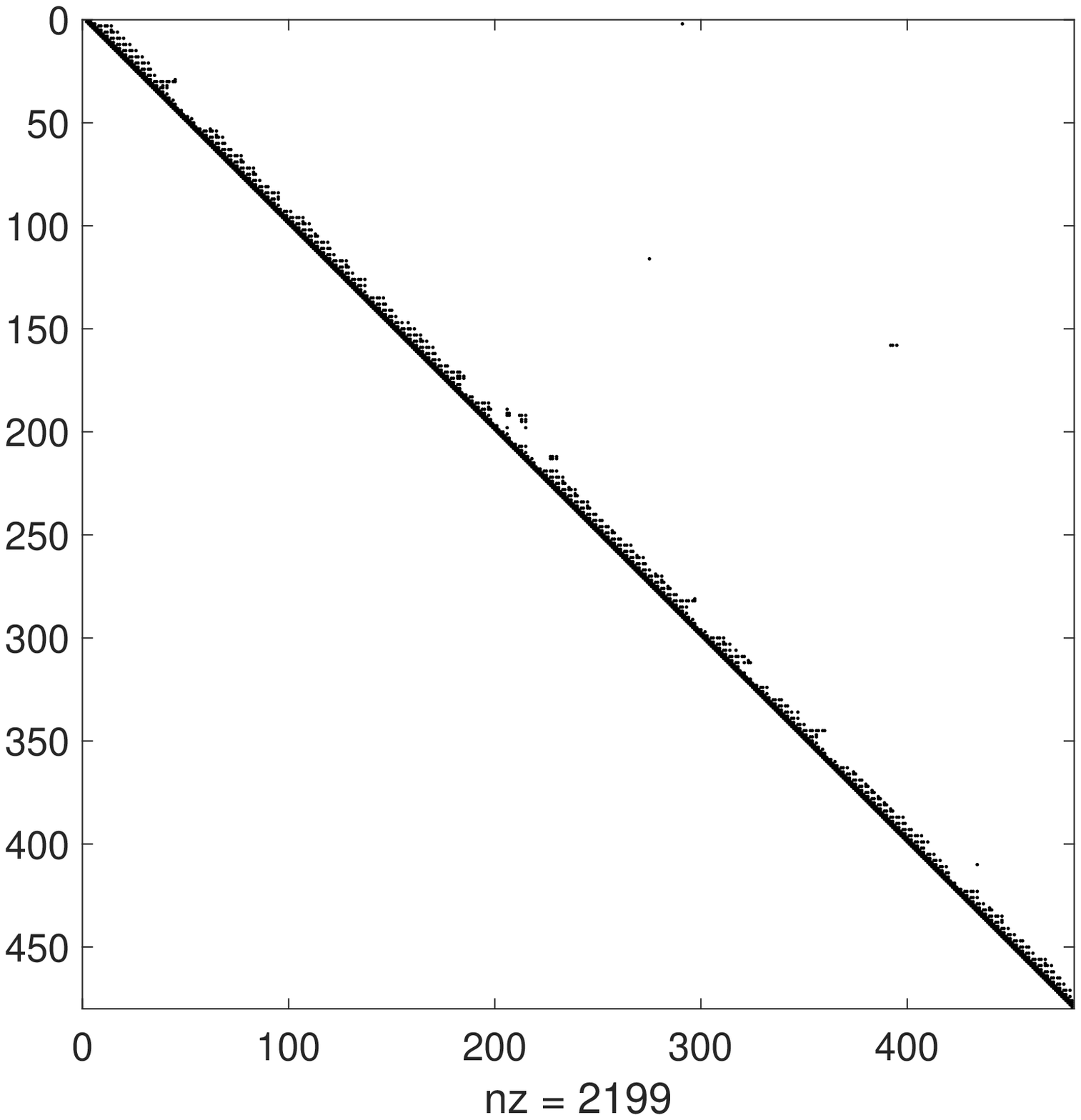}
\caption{Two typical distance distributions for protein-like instances: 1B4R (left) and 1ALL (right).}\label{fig:dist}
\end{figure}

The natural backbone order for instances generated in this way provides a vertex order satisfying the assumptions of Definition~\ref{def:dmdgp}, 
implying we are working with $^3$DMDGP instances. 
In Figure~\ref{fig:dist}, we show the known entries of distance matrices for this kind of problems.

In our experiments we consider two test sets: one using cut-off 6 \AA \ and other using 5 \AA \ . 
In Tables~\ref{tab:results6} and \ref{tab:results5}, we present the results obtained by both algorithms: 
%BP is the classical Branch-and-Prune implementing a depth-first search \cite{llm08,coap2012}, 
%where no symmetry is exploited at all and coordinates 
%for the vertices are recomputed every time a backtracking occurs, 
%whereas  Algorithm~\ref{algo:SBBU2}, called SBBU (Symmetry-based Build-up), 
%uses a pruning edge based search (described in Section~\ref{sec:algo}), 
%along with merging and partial reflections of sub-structures to properly exploit symmetries 
%of DMDGP instances to find the first solution. 
BP is the classic Branch-and-Prune implementing a depth-first search \cite{llm08,coap2012}, 
%where no symmetry is exploited at all and coordinates 
%for the vertices are recomputed every time a backtracking occurs, 
whereas  SBBU (Symmetry-based Build-up) corresponds to Algorithm~\ref{algo:SBBU2}.
%uses a pruning edge based search (described in Section~\ref{sec:algo}), 
%along with merging and partial reflections of sub-structures to properly exploit symmetries 
%of DMDGP instances to find the first solution. 

The algorithms were implemented in C++ and the experiments carried out 
in Intel(R) Core(TM) i7-3770 CPU @ 3.40GHz, 8G RAM, 
running Linux Ubuntu 18.04.4, gcc version 7.4.0 compiler. 

%\tred{TODO: include comments and links to GitHub repo provided by Antonio? Will we make our code also public available?}

Both tables bring the ID of the protein in PDB, the number of atoms $|V|$, number of edges (available distances) $|E|$, 
CPU time in seconds for the two algorithms and the normalized Mean Distance Deviation (MDE): 
\begin{equation}\label{eq:mde}
MDE(X,E,d) = \dfrac{1}{|E|} \sum_{\{i,j\} \in E} \dfrac{\left| \, \| x_i - x_j \|_2 - d_{ij} \, \right|}{d_{ij}}.
\end{equation}
Both algorithms were stopped as soon as the first solution is found and 
a ``--'' symbol means that the algorithm was not able to find a solution in less than 300 seconds. 
For each instance, we also present the total and maximum works $W$ and $\bar{W}$, respectively.   
%for the specific pruning edge order we are using in this paper. 
The last column, called ``Speed-up'', contains the ratio time BP / time SBBU. 

From the tables, we observe that the new algorithm provides a non-trivial speed-up 
in most of the instances. In particular, the new algorithm is considerably faster than the classic BP 
for the sparser instances where it was up to 1,000 times faster.

\begin{table}[]
\setlength{\tabcolsep}{5pt}
\footnotesize
\centering 
\begin{tabular}{|lrr|rr|rrrr|r|}
\hline 
 &  &  & \multicolumn{2}{c|}{BP} &  \multicolumn{4}{c|}{SBBU} &  \\ \hline 
ID & $|V|$ & $|E|$ & Time & MDE & Time & MDE & $\bar{W}$ & $W$ & Speed-up \\ \hline 
1N6T & 30 & 236 & 7.60E-05 & 8.32E-05 & 1.77E-05 & 2.72E-12 & 2 & 52 & 4.29 \\
1FW5 & 60 & 558 & 1.30E-04 & 1.51E-05 & 3.51E-05 & 4.22E-12 & 2 & 112 & 3.70 \\
1ADX & 120 & 1008 & 2.10E-04 & 5.62E-12 & 4.49E-05 & 3.78E-12 & 2 & 232 & 4.67 \\
1BDO & 241 & 2167 & 4.10E-04 & 3.79E-12 & 9.24E-05 & 1.39E-11 & 2 & 474 & 4.44 \\
1ALL & 480 & 4932 & 8.40E-04 & 8.91E-13 & 1.90E-04 & 3.80E-12 & 2 & 952 & 4.42 \\
6S61 & 522 & 5298 & 8.70E-04 & 6.50E-13 & 2.06E-04 & 3.09E-12 & 2 & 1036 & 4.23 \\
1FHL & 1002 & 9811 & 2.00E-03 & 6.82E-12 & 3.97E-04 & 1.93E-11 & 2 & 1996 & 5.04 \\
4WUA & 1033 & 9727 & 1.80E-03 & 1.47E-11 & 3.94E-04 & 7.73E-12 & 8 & 2060 & 4.57 \\
6CZF & 1494 & 14163 & 2.60E-03 & 1.33E-12 & 5.79E-04 & 4.18E-12 & 2 & 2980 & 4.49 \\
5IJN & 1950 & 18266 & 3.40E-03 & 1.37E-12 & 7.64E-04 & 1.76E-11 & 16 & 3908 & 4.45 \\
6RN2 & 2052 & 19919 & 3.70E-03 & 1.11E-12 & 8.27E-04 & 1.54E-11 & 16 & 4104 & 4.48 \\
1CZA & 2694 & 26452 & 4.90E-03 & 1.29E-12 & 1.07E-03 & 6.22E-11 & 2 & 5380 & 4.59 \\
6BCO & 2856 & 27090 & 7.90E-03 & 4.53E-13 & 1.10E-03 & 7.91E-12 & 16 & 5730 & 7.15 \\
1EPW & 3861 & 35028 & 7.80E-03 & 1.88E-11 & 1.44E-03 & 2.50E-10 & 2 & 7714 & 5.40 \\
5NP0 & 7584 & 80337 & 3.10E-02 & 6.58E-12 & 3.58E-03 & 1.35E-10 & 256 & 15562 & 8.66 \\
5NUG & 8760 & 82717 & 2.40E-02 & 1.43E-06 & 3.45E-03 & 5.33E-10 & 16 & 17592 & 6.96 \\
4RH7 & 9015 & 85831 & 2.50E-02 & 1.62E-12 & 3.67E-03 & 2.22E-10 & 16 & 18054 & 6.82 \\
3VKH & 9126 & 87621 & 2.70E+00 & 3.00E-08 & 3.62E-03 & 1.15E-09 & 256 & 18556 & 745.03 \\ \hline 

\end{tabular}
% resultados com cutoff 6A
\caption{Computational results in some protein-like instances (cut-off: 6\AA ).}\label{tab:results6}
\end{table}

\begin{table}[]
\setlength{\tabcolsep}{5pt}
\footnotesize
\centering 
\begin{tabular}{|lrr|rr|rrrr|r|}
\hline 
 &  &  & \multicolumn{2}{c|}{BP} &  \multicolumn{4}{c|}{SBBU} &  \\ \hline 
ID & $|V|$ & $|E|$ & Time & MDE & Time & MDE & $\bar{W}$ & $W$ & Speed-up \\ \hline 
1N6T & 30 & 176 & 7.60E-05 & 5.14E-05 & 1.04E-05 & 5.64E-12 & 2 & 52 & 7.31 \\ 
1FW5 & 60 & 417 & 1.40E-04 & 7.99E-06 & 2.11E-05 & 3.08E-12 & 2 & 112 & 6.63 \\
1ADX & 120 & 659 & 4.70E-04 & 3.50E-06 & 3.49E-05 & 2.53E-12 & 2 & 232 & 13.48 \\
1BDO & 241 & 1345 & 3.60E-04 & 1.50E-07 & 7.05E-05 & 1.04E-11 & 2 & 474 & 5.11 \\
1ALL & 480 & 3443 & 9.80E-04 & 2.81E-06 & 1.67E-04 & 1.27E-12 & 2 & 952 & 5.88 \\
6S61 & 522 & 3699 & 8.70E-04 & 8.10E-07 & 1.75E-04 & 1.39E-12 & 2 & 1036 & 4.98 \\
1FHL & 1002 & 6378 & 2.70E-03 & 2.56E-12 & 2.88E-04 & 1.17E-11 & 2 & 1996 & 9.38 \\
4WUA & 1033 & 6506 & 1.80E-03 & 5.34E-12 & 2.96E-04 & 2.94E-12 & 16 & 2066 & 6.08 \\
6CZF & 1494 & 9223 & 2.40E-03 & 4.62E-13 & 4.36E-04 & 2.33E-12 & 2 & 2980 & 5.51 \\
5IJN & 1950 & 11981 & 4.00E-03 & 4.43E-13 & 6.08E-04 & 4.23E-12 & 16 & 3908 & 6.58 \\
6RN2 & 2052 & 13710 & 5.50E-03 & 3.89E-13 & 8.58E-04 & 9.35E-12 & 16 & 4112 & 6.41 \\
1CZA & 2694 & 17451 & 5.80E-03 & 4.51E-13 & 8.03E-04 & 3.06E-11 & 2 & 5380 & 7.22 \\
6BCO & 2856 & 18604 & 5.00E-03 & 6.00E-07 & 1.05E-03 & 6.96E-12 & 16 & 5706 & 4.75 \\
1EPW & 3861 & 23191 & 2.30E-02 & 3.00E-08 & 1.13E-03 & 9.78E-11 & 8 & 7716 & 20.29 \\
5NP0 & 7584 & 59478 & 2.90E-01 & 2.56E-12 & 2.80E-03 & 4.11E-11 & 256 & 16138 & 103.55 \\
5NUG & 8760 & 56979 & 2.70E+00 & 3.60E-07 & 2.67E-03 & 1.05E-10 & 128 & 17700 & 1011.09 \\
4RH7 & 9015 & 59346 & 3.10E-02 & 5.64E-13 & 2.97E-03 & 1.20E-10 & 32 & 18068 & 10.43 \\
3VKH & 9126 & 59592 & -- & -- & 2.45E-02 & 1.10E-09 & 65536 & 84066 & \\ \hline 
\end{tabular}
% resultados com cutoff 5A
\caption{Computational results in some protein-like instances (cut-off: 5\AA ).}\label{tab:results5}
\end{table}

\begin{figure}
\centering
\includegraphics[scale=0.18,trim=30 0 30 0,clip]{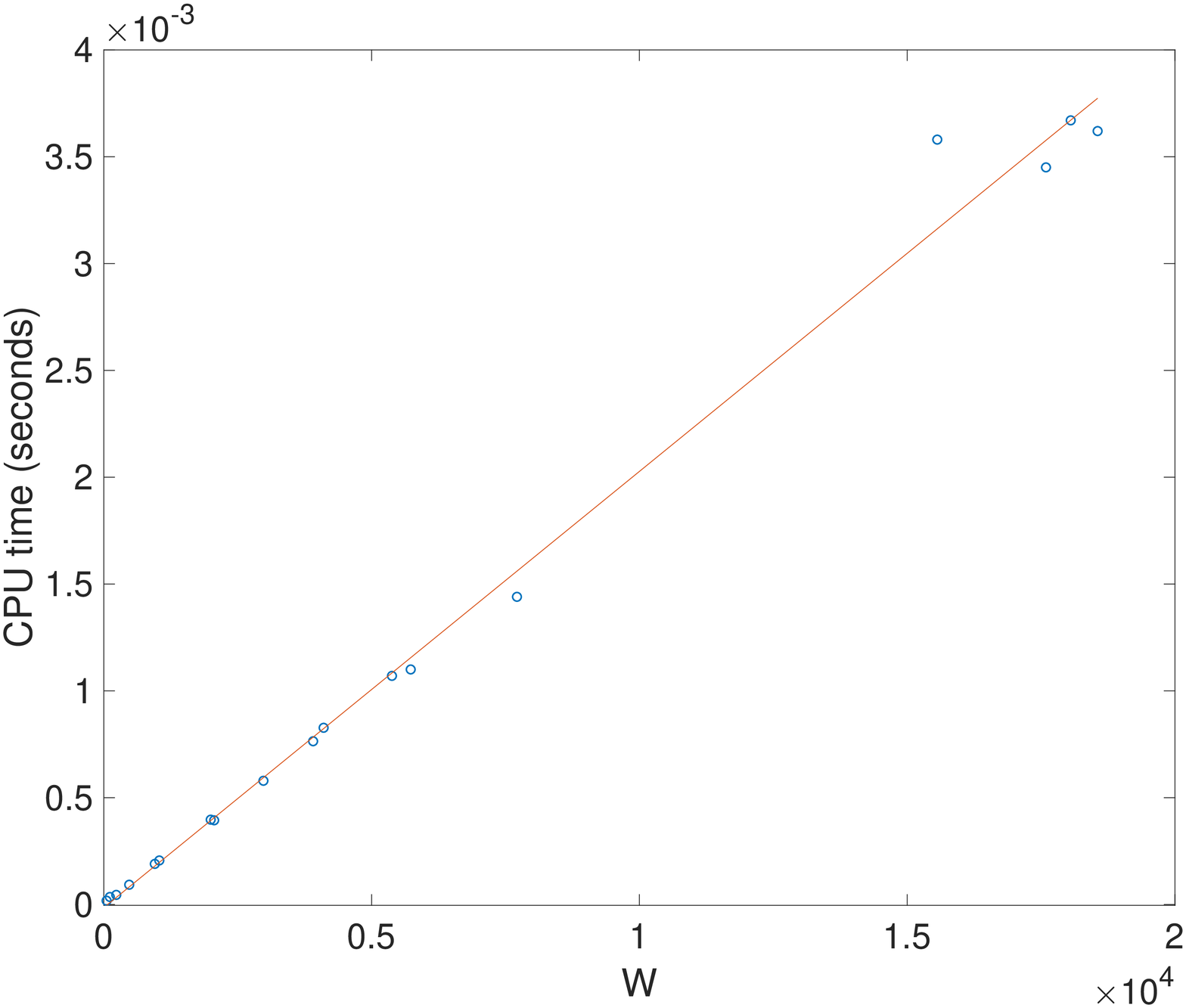}
\includegraphics[scale=0.18,trim=30 0 30 0,clip]{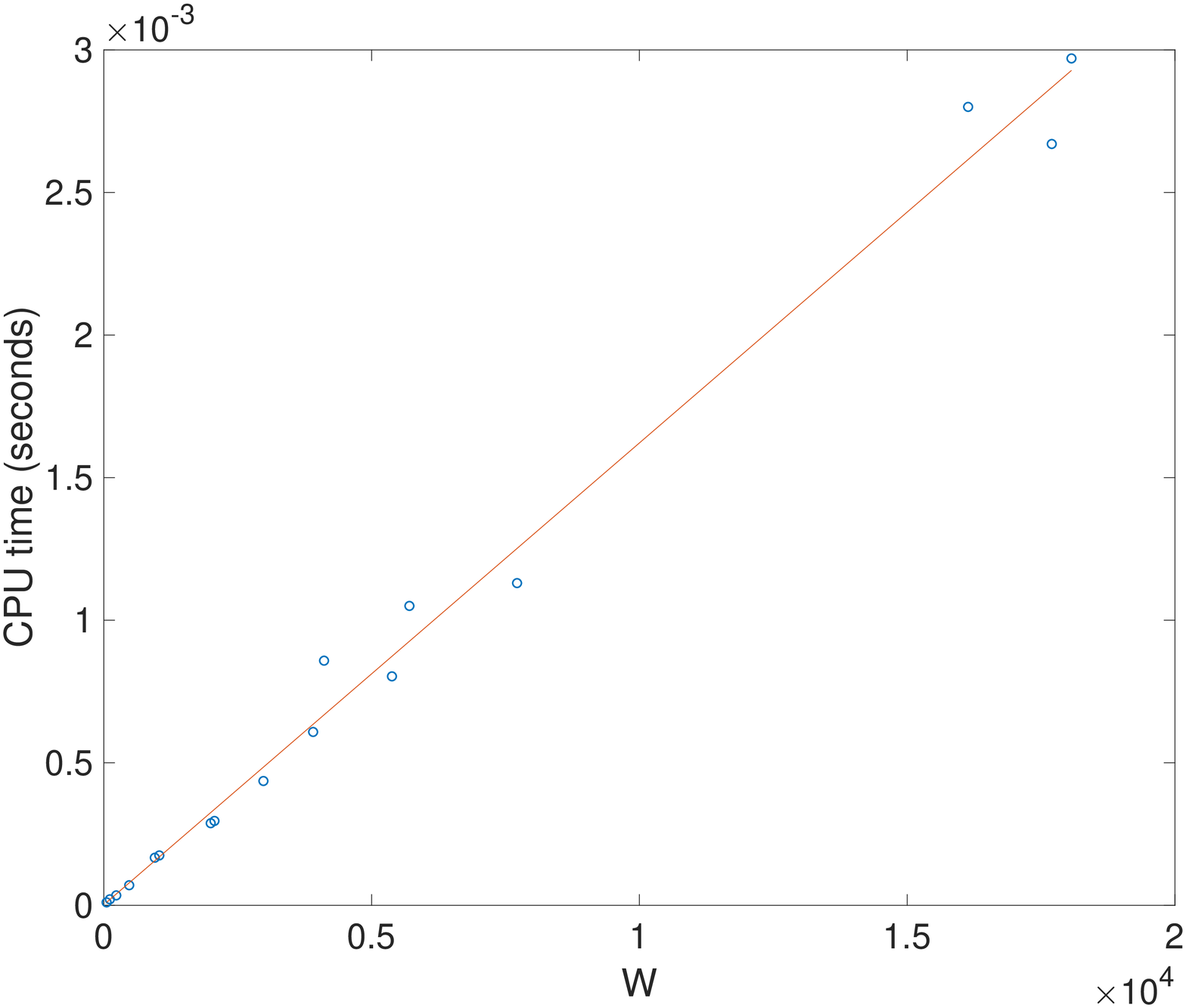}
\caption{Scatter plots $W \times$ time and linear regressions for the two datasets 
(Table~\ref{tab:results6} on the left, Table~\ref{tab:results5} on the right but not considering the last row).}\label{fig:linear}
\end{figure}

Concerning the estimated total work of SBBU, it seems that the time varies linearly with $W$ 
as depicted in Figure~\ref{fig:linear}. The relationship between BP time and $W$ and/or $\bar{W}$ 
is not so clear. However, we argue that while the most costly subproblem $(G^{ij},K)$ represents a cost of $\bar{W}$ 
in the total cost $W$ for SBBU, for the usual DFS recursive implementation of BP, it may contribute 
much more to the BP total cost because such subproblem may have to be solved several times 
in the occasion of backtrackings. 

\section{Conclusion}\label{sec:final}

We propose a new algorithm for the DMDGP which leverages symmetry information to find the first solution quickly. 
By efficiently exploiting symmetries of subproblems defined by pruning edges, 
and reducing the degrees of freedom by taking into account already solved subproblems, 
%merging already solved subproblems in rigid bodies, 
%and following a specific edge order, 
the resulting algorithm  appears to be quite efficient in sparse DMDGP instances, sometimes giving a significant speed-up with respect to the classic BP algorithm. 

In the proposed version of SBBU algorithm we consider a specific order for the pruning edges. 
In future works we expect to generalize SBBU in order to handle different pruning edges orderings 
and study the impact of these in the total cost $W$. 

\section{Acknowledgments}

%We would like to thank Prof. Luiz M. Carvalho for valuable discussions and
%to the Brazilian research agencies CNPq, CAPES, and FAPESP for the finantial
%support. Part of this work was done during the visit of DG to LL 
%at École Polytechnique, supported by CAPES/Print Process 88881.310538/2018-01.  

We would like to thank Prof.~Luiz M.~Carvalho for valuable discussions. 
This work was partly supported by: 
(a) the Brazilian research agencies CNPq, CAPES, and FAPESP; 
(b) the French research agency ANR under grant ANR-19-CE45-0019 ``multiBioStruct"; 
(c) the European Union's Horizon 2020 research and innovation programme under the Marie Sklodowska-Curie grant agreement n. 764759 ETN "MINOA". 
Part of this work was done during the visit of DG to LL at \'{E}cole Polytechnique, supported by CAPES/Print Process 88887.465828/2019-00.

\bibliographystyle{plain}
\bibliography{references}

\section*{Appendix}
Proof of Lemma~\ref{lem:NaoEstraga}. 
\begin{proof}
Since $x(s')$ from Eq.~\eqref{eq:xs} involves only partial reflections, 
in view of Property 3 in Remark~\ref{rem:pr}, $x(s') \in \hat{X}$, i.e~$\| x_u(s') - x_w(s') \| = d_{uw}, \forall \{u,w\} \in E_D$. 

It remains to show that $x(s')$ does not violate distance constraints associated to pruning edges in $P^{ij}$. 
Since the reflections are applied to positions $x_{\ell}$ such that $\ell \ge i+K+1$, 
pruning edges $\{u,w\} \in P^{ij}$ with $u<w\le i+K$ are not affected. 
Thus, assume that $i+K+1 \le w \le n$. 
If $u \le i$, then for $\ell = i+K+1,\dots,w$ there exists $\{u,w\}$ 
such that $u+K+1 \le \ell \le w$, which implies that $v_{i+K+1},\dots,v_{w} \not\in S^{ij}$, 
meaning that the first symmetry vertex $v_{\ell}$ in $S^{ij}$ is such that $\ell \ge w+1$. 
Thus, according to \eqref{eq:g} and \eqref{eq:xs}, partial reflections are not applied to $x_{i+K+1},\dots,x_{w}$, 
i.e~$x_{\ell}(s') = x_{\ell}(s)$, for $\ell=i+K+1,\dots,w$ and $\| x_u(s') - x_w(s') \| = d_{uw}$ holds. 
Otherwise, for $u \ge i+1$, we have that $v_{u+K+1},\dots,v_w \not\in S^{ij}$, and from \eqref{eq:xs} and \eqref{eq:g}, 
positions $x_{\ell},\dots,x_{u+K+1},\dots,x_w$ are updated by reflections 
$R_{x}^{\ell}(x_{\ell}), \dots, R_{x}^{\ell}(x_{u+K+1}), \dots, R_{x}^{\ell}(x_w)$, for $i+K+1 \le \ell \le u+K$ such that $v_{\ell} \in S^{ij}$. 
Since either $u \le \ell - 1$, i.e~$x_u$ is in the hyperplane associated to $v_{\ell}$, 
or $u \ge \ell$, i.e.~$x_u$ comes after this hyperplane, 
in view of Remark~\ref{rem:pr}, Property 2, these reflections are such that $\| x_u(s') - x_v(s') \| = d_{uw}$. 
\end{proof}

\end{document}